\documentclass[pdflatex,sn-mathphys-num]{sn-jnl}

\usepackage{graphicx}%
\usepackage{multirow}%
\usepackage{amsmath,amssymb,amsfonts}%
\usepackage{amsthm}%
\usepackage{mathrsfs}%
\usepackage[title]{appendix}%
\usepackage{textcomp}%
\usepackage{manyfoot}%
\usepackage{booktabs}%
\usepackage{algorithm}%
\usepackage{algorithmicx}%
\usepackage{algpseudocode}%
\usepackage{listings}%

\usepackage{physics}
\usepackage{mathtools}
\usepackage{bm}
\mathtoolsset{showonlyrefs}
\DeclareMathOperator{\supp}{supp}
\DeclareMathOperator{\diag}{diag}

\theoremstyle{thmstyleone}%
\theoremstyle{thmstyletwo}%
\newtheorem{remark}{Remark}%

\theoremstyle{thmstylethree}%

\begin{document}

\title[Article Title]{A fast direct solver for two-dimensional transmission problems of elastic waves}


\author*[1]{\fnm{Yasuhiro} \sur{Matsumoto}}\email{matsumoto@cii.isct.ac.jp}

\author[2]{\fnm{Taizo} \sur{Maruyama}}\email{maruyama.t.45ef@m.isct.ac.jp}


\affil*[1]{\orgdiv{Center for Information Infrastructure}, \orgname{Institute of Science Tokyo}, \orgaddress{\street{2-12-1-I8-21, Ookayama}, \city{Meguro-ku}, \postcode{152-8550}, \state{Tokyo}, \country{Japan}}}

\affil[2]{\orgdiv{Department of Civil and Environmental Engineering}, \orgname{Institute of Science Tokyo}, \orgaddress{\street{2-12-1-W8-22, Ookayama}, \city{Meguro-ku}, \postcode{152-8550}, \state{Tokyo}, \country{Japan}}}



\abstract{
This paper describes a fast direct boundary element method for elastodynamic transmission problems in two dimensions,
which can be used for analyzing elastic wave scattering by an inclusion.
We develop an efficient solver based on a discretization method that is
broadly applicable regardless of the inclusion shape.
From the smoothness of the solutions of the Navier--Cauchy equation,
it is reasonable that
the displacement is approximated by the piecewise linear bases
and the traction is approximated by the piecewise constant bases.
However, in this mixed bases strategy, Calder\'on preconditioning, that
is, an analytical preconditioning with excellent performance, cannot be applied.
To circumvent this issue,
we developed a fast direct solver formulated using both Burton--Miller and Poggio--Miller--Chang--Harrington--Wu--Tsai (PMCHWT) boundary integral equations.
Our method uses a technique based on the proxy method for low-rank approximation of the coefficient matrix's off-diagonal blocks. 
To handle transmission problems, the proposed fast direct solver uses separate binary tree partitions for nodes and elements.
Numerical examples demonstrate that our solver achieves linear computational complexity at fixed low frequencies
and 
can efficiently handle problems with multiple right-hand sides.
Notably, the solver based on the Burton--Miller formulation is approximately 20\% faster than the one using the PMCHWT formulation.
Our new method provides a versatile, fast solver, whose performance is relatively independent of the shape of inclusions and computational parameters, such as density, for elastodynamic transmission problems.
}




\maketitle

\section{Introduction}\label{sec:intro}
The elastic wave scattering problem is important for
mechanical engineering, architecture, and civil engineering,
because it often arises in ultrasonic non-destructive testing \cite{CHEN2021123832},
seismic wave analysis \cite{bouchon2007boundary, antonietti2018numerical}, and environmental vibration evaluation \cite{wang2025review, Avillez2021Procedures}.
The domains dealt with are not always composed of a single material.
Recently,
composite materials like those used in aircraft are frequently adopted for their light weight and strong structures.
Concrete, a long-used material, is also a composite because the wave speed and density in the cement paste and aggregate are different.
The propagation of elastic waves in composite materials can be formulated as an elastodynamic transmission problem.
This problem generally needs to be solved numerically because analytical solutions can be obtained only for some simple settings.
Inclusions that are relatively large in comparison to the wavelength lead to an increased number of degrees of freedom (DOFs) in discrete systems.
Therefore, we need an efficient numerical method to solve large-scale transmission problems.

This study considered time-harmonic wave scattering.
For linear scattering problems, time-domain solutions can be constructed by time-harmonic solutions using the inverse Fourier transform.
Numerical methods based on the boundary integral equations,
referred to as boundary element methods in this work,
are suitable for time-harmonic wave scattering \cite{kobayashi_book}.
In contrast to approaches using the finite element method or finite difference method,
boundary element methods can deal with the radiation condition
at infinity without requiring specific handling \cite{liu2012recent}.
A drawback of the boundary element method is that the coefficient matrix for its discretized system is dense.
The dense matrix leads to large computational costs for its construction and solving the system of linear equations.
This drawback is accentuated in large-scale problems.

The problem with the computational cost for dense systems in the boundary element methods is well known, and several solutions have been proposed.
One solution uses high-order discretization methods, such as isogeometric analyses \cite{nguyen2015isogeometric} or Nystr\"om methods \cite{canino1998numerical}.
Boundary element methods formulated in the isogeometric analysis (e.g., elasticity \cite{simpson2012two} and acoustics \cite{simpson2014acoustic})
represent the surfaces of inclusions of a computational region as a linear combination of high-order polynomials.
The same high-order polynomials are also employed for the basis of the approximated unknown functions.
High-order quadrature formulas for polynomials representing the surface have been proposed to handle weak-, strong-, and hyper-singular integrands
that are calculated in boundary integral equations \cite{calabro2018efficient, aimi2020quadrature}.
A combination of the isogeometric analysis and the hierarchical matrix method \cite{borm2003introduction} has also been proposed \cite{desiderio2025hierarchical}.
However, different implementations are needed for these methods,
depending on the kind of boundary geometry (e.g., smooth boundary or boundary with corner points).
Furthermore, the elastodynamic versions of quadrature formulas for the isogeometric analysis are not easy to construct because their integrand functions are complicated.

Another choice for high-order discretization methods is Nystr\"om methods for elastodynamics \cite{TONG20071845}.
It is a collocation method, whose collocation points and integration points coincide, using high-order quadrature formulas for singular integral kernels.
As mentioned above, because constructing a high-order quadrature formula for elastodynamics is difficult,
there is no definitive method, to the best of the authors' knowledge.
Recently, Nystr\"om methods based on the Helmholtz decomposition \cite{guo2025quadrature}
and the quadrature by expansion \cite{dominguez2024nystrom}
were proposed.
However, with the Nystr\"om method, just as with isogeometric boundary element methods,
the quadrature method needs to be changed depending on the boundary geometry.
For example, boundaries with corner points necessitate the use of sigmoid transformations or graded mesh quadratures
tailored to those corners \cite{kress1990nystrom, dominguez2016well, bremer2012nystrom},
whereas such treatments are not required for smooth boundaries.

From a practical perspective, it is important that a program code can be used for general conditions.
This work aimed to develop an efficient solver for elastodynamic transmission problems based on a discretization method broadly applicable to a computational region that has inclusions with a closed Lipschitz boundary.
One of the easiest way to deal with Lipschitz domains is to use a low-order Galerkin discretization,
because we can approximate the boundary of an inclusion with corners as a polygon.
Considering the continuity of displacement and discontinuity of stress on the surface,
it is reasonable to approximate the displacement with piecewise linear bases
and approximate the traction with piecewise constant bases on this polygon.
However, with this Galerkin discretization method using both piecewise linear and constant bases as test and basis functions,
it is non-trivial how to apply the Calder\'on preconditioning, which is a famous analytical preconditioning method,
for an iterative solver such as GMRES \cite{niino2018calderon}.
Moreover, it is difficult to take advantage of the fast multipole method (FMM) \cite{yoshida2001application},
which is a typical fast method used in iterative solvers.
We note that if the boundary is smooth,
  the Calder\'on preconditioning on the Galerkin discretization using only the piecewise linear basis can be applied with an additional preconditioner formed as the Gram matrix associated with basis functions \cite{NIINO201266, Wout2021Benchmarking}.
Besides Calder\'on preconditioning, there are also several other analytical preconditioning methods that can potentially outperform it \cite{Wout2021Benchmarking}.
For instance, on-surface radiation condition (OSRC) methods \cite{CHAILLAT2017429} achieve preconditioning of boundary integral equations by utilizing the Dirichlet-to-Neumann or Neumann-to-Dirichlet maps.
Although there is an example of the OSRC method that the boundary integral equations for the Helmholtz transmission problem in a Lipschitz domain were efficiently solved \cite{VANTWOUT2022111229}, this approach approximated both the solution and its normal derivative using only the piecewise linear basis function.
Furthermore, we note that transferring these successful preconditioning methods from Helmholtz scattering problems to elastic wave scattering problems may not be straightforward, due to the difference in the singularity of the double-layer kernel \cite{BRUNO2020109350}.
When analytical preconditioning methods are not applicable, we can use numerical preconditioning techniques based on
methods utilizing FMM itself for preconditioning \cite{Carpentieri2005Combining},
algebraic multigrid methods \cite{Di2023Algebraic},
and structured hierarchical matrices \cite{Spendlhofer2025ISC}.
However it is documented that analytical preconditioners often outperform algebraic preconditioners \cite{Antoine01102008, VANTWOUT2022111229}.

As an alternative for solving the problem under this discretization scheme, fast direct solvers are appealing.
Some of the shortcomings of the iterative solvers can be overcome by utilizing fast direct solvers \cite{martinsson2019book}.
In this work, we limited our scope to two-dimensional problems.
This is because, while the discretization methods discussed in this work are applicable in both two and three dimensions,
constructing a fast direct solver is significantly harder in three dimensions than in two dimensions.
Fast direct solvers for solving three-dimensional problems have been proposed \cite{Greengard2009fast, chaillat2017theory, rong2019fast, sushnikova2023fmm, ma2024inherently}, but they are still under development.
The primary challenge lies in approximating the off-diagonal blocks of the discrete system's coefficient matrix to a structured low-rank form and integrating this into a fast direct solver.
In three dimensions, strong admissibility \cite{minden2017recursive} for low-rank approximation is required,
and its incorporation into fast direct solvers remains an active area of research.
Fortunately, this issue is mitigated in two dimensions, as low-rank approximation based on weak admissibility works effectively.
Weak admissibility means that, when we evaluate interactions between cells in a binary tree structure introduced into a computational region, all interactions between different cells are approximated as low rank.
Conversely, strong admissibility approximates only far-field interactions between different cells to low rank, leaving near-field interactions dense.

This work developed a fast direct solver, which is an extension of the case of elastic wave scattering by a cavity \cite{MATSUMOTO2025106148} to the case of transmission problems with one homogeneous inclusion.
It is a variant of the Martinsson--Rokhlin fast direct solver based on the proxy method \cite{MARTINSSON20051}
and achieves $O(N)$ complexity at fixed low frequencies in two dimensions, where $N$ is the DOFs of the discretized system.
In elastodynamics, few studies have addressed fast direct solvers, to the best of our knowledge.
There is the application of the standard LU factorization
of $\mathcal{H}$-matrices \cite{chaillat2017theory},
which has $O(N \log^2 N)$ complexity in the best case.
We formulated a fast direct solver
under the discretization method using both piecewise linear and constant test/basis functions.
To handle transmission problems, the proposed fast direct solver
uses separate tree partitions for nodes and elements.
We consider not only the Poggio--Miller--Chang--Harrington--Wu--Tsai (PMCHWT) formulation \cite{POGGIO1973159, Chang1977surface}, which is the usual choice,
but also the formulation based on the Burton--Miller method \cite{burton1971application}.
For the considered two subdomains case, the Burton--Muller formulation has six layer potentials, while the PMCHWT formulation has eight layer potentials.
Therefore, the Burton--Muller formulation is expected to be faster than the PMCHWT formulation.
This paper explains the weak admissibility type low-rank approximation method using the proxy method for both formulations in detail.
The program code for the proposed method is versatile because it can be used whether a geometry includes or does not include corner points.
Other options for the boundary integral equation for transmission problems include
the single-boundary integral equation \cite{kleinman1988single} and the multi-trace-boundary integral equation \cite{hiptmair2012multiple, Claeys2013multi}.
However, neither are suitable for the proposed fast direct solver because the former uses an indirect unknown density and products of layer potentials, and the number of unknowns in the latter method is twice that in the PMCHWT or Burton--Miller formulation, although it uses direct (physical) unknowns.
We use numerical examples to show that the proposed fast direct solver has $O(N)$ complexity at fixed low frequencies and is relatively robust with respect to the computational time for changing calculation conditions such as the density of an inclusion.
Furthermore, the proposed fast direct solver based on the Burton--Miller method outperforms that based on the PMCHWT formulation by about $20\%$ with respect to the computational time.

The remainder of this manuscript is organized as follows.
A statement of the transmission problems of elastic waves and corresponding boundary integral equations is given in
Section \ref{sec:statement}.
Section \ref{sec:Discretization} discusses the Galerkin discretization of the boundary integral equations and its block formulation.
Section \ref{sec:proposed} discusses the proposed method comprising the proxy method and the compression technique of the linear equations.
Section \ref{sec:Numerical} demonstrates the correctness and performance of the proposed fast direct solver by numerical examples.
Section \ref{sec:Conclusion} concludes this work and mentions future work.

\section{Statement of problems and boundary integral equations} \label{sec:statement}
\subsection{Transmission problems for elastic waves}
\begin{figure}[tb]
  \centering
  \includegraphics[width=0.35\linewidth]{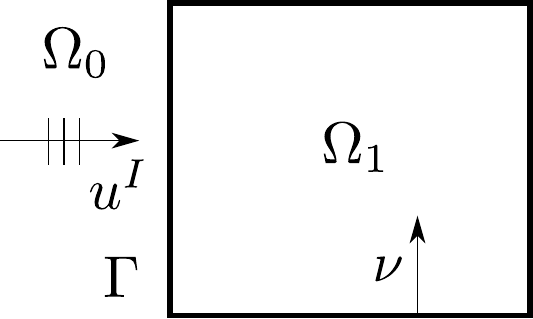}
  \caption{Diagram of transmission problems for an elastic wave}
  \label{fig:domain}
\end{figure}
Unless otherwise specified, the Cartesian coordinate system $x = (x_1, x_2) \in \mathbb{R}^2$ is used to express the wave field.
We consider a time-harmonic in-plane elastic wave scattering problem
with time factor $\exp(-\mathrm{i} \omega t)$,
where $\mathrm{i}$, $\omega$, and $t$ are the imaginary unit, angular frequency, and time, respectively.
Consider a bounded domain $\Omega_1 \subset \mathbb{R}^2$ whose boundary $\Gamma = \partial \Omega_1$
is a Jordan closed curve and is Lipschitz.
We interpret $\Omega_1$ as an inclusion.
Let $\Omega_0 = \mathbb{R}^2 \setminus \overline{\Omega_{1}}$ be a (connected) domain filled by a linearly elastic, homogeneous, and isotropic solid with density $\rho^{(0)}$ and Lam\'e constants $\lambda^{(0)}$ and $\mu^{(0)}$.
Assume that $\Omega_1$ is also filled by a linearly elastic, homogeneous, and isotropic solid with density $\rho^{(1)}$ and Lam\'e constants $\lambda^{(1)}$ and $\mu^{(1)}$.
In $\Omega_0$, there is an incident wave
 $u^{I}(x) := (u_{1}^{I}(x), u_{2}^{I}(x)) : \mathbb{R}^2 \to \mathbb{C}^2$,
which is a vector field.
Let $u(x) := (u_1 (x), u_2 (x)) : \mathbb{R}^2 \to \mathbb{C}^2$ be the total
displacement vector field and $u^{S}_{i}(x) := u_{i}(x) - u^{I}_{i}(x)$
be the scattered wave field in $\Omega_{0}$ and $\Omega_{1}$.  Note that
$u(x) \in \Gamma$ is defined by \eqref{eq:transmission_condition}.
Let $\nu(x) = (\nu_{1}(x), \nu_{2}(x)) : \mathbb{R}^2 \to \mathbb{R}^2$ be a unit normal vector
on $\Gamma$ toward $\Omega_1$.  The above setting is shown
 in Figure \ref{fig:domain}.
For a function $v(x)$, let the notations $v_{,i}(x)$ and $v_{,ij}(x)$ be defined as
\begin{equation}
  v_{,i}(x) := \pdv{v(x)}{x_i}, \quad v_{,ij}(x) := \pdv{v(x)}{x_i}{x_j}.
\end{equation}

The transmission problem for elastic waves is to find the solution $u(x)$ that satisfies the
Navier--Cauchy equation in $\Omega_{0}$ and $\Omega_{1}$:
\begin{multline}
  \mu^{(m)} \sum_{j=1}^{2} u_{i,jj}(x) + (\lambda^{(m)} + \mu^{(m)}) \sum_{j=1}^{2} u_{j, ji}(x) + \rho^{(m)} \omega^{2} u_{i}(x) = 0,\\ x \in \Omega_{m}, \quad m = 0, 1, \quad i = 1, 2, \label{eq:navie}
\end{multline}
with the transmission boundary conditions
\begin{equation}
  \begin{dcases}
    \lim_{h \downarrow 0} u_{i}(x + h \nu(x)) = \lim_{h \downarrow 0} u_{i}(x - h \nu(x)) \quad \qty( =: u_i (x)), \quad x \in \Gamma, \quad i = 1, 2, \\
    \begin{split}
    \lim_{h \downarrow 0}\sum_{j, p, q=1}^{2} C_{ipjq}^{(1)} u_{j, q}(x + h \nu(x)) \nu_{p}(x) = \lim_{h \downarrow 0}\sum_{j, p, q=1}^{2} C_{ipjq}^{(0)} u_{j, q}(x - h \nu(x)) \nu_{p}(x) \\
   \qty( =: t_i (x)), \quad x \in \Gamma, \quad i = 1, 2,
    \end{split}
  \end{dcases}
  \label{eq:transmission_condition}
\end{equation}
and under the radiation conditions and the regularity conditions of elastodynamics for $u^{S}(x)$ \cite{eringen1975elastodynamics}:
\begin{equation}
  \begin{dcases}
    \lim_{|x| \to \infty} |x|^{\frac{1}{2}} \qty( \pdv{u^{S}_{i;L}(x)}{|x|} - \mathrm{i}k_{L}^{(0)} u^{S}_{i;L} (x)) = 0, \\
    \lim_{|x| \to \infty} |x|^{\frac{1}{2}} \qty( \pdv{u^{S}_{i;T}(x)}{|x|} - \mathrm{i}k_{T}^{(0)} u^{S}_{i;T} (x)) = 0, \\
    \lim_{|x| \to \infty} |x|^{-\frac{1}{2}} u^{S}_{i;L} (x) = 0, \\
    \lim_{|x| \to \infty} |x|^{-\frac{1}{2}} u^{S}_{i;T} (x) = 0.
  \end{dcases}
  \quad (\text{for } i = 1, 2)
  \label{eq:raditation}
\end{equation}
Here, $t_{i}(x)$ is the traction in the $x_{i}$ direction;
$C_{ipjq}^{(m)} := \lambda^{(m)} \delta_{ip}\delta_{jq} + \mu^{(m)} (\delta_{ij}\delta_{pq} + \delta_{iq}\delta_{pj})$ is a component of the elastic tensor; $\delta_{ij}$ is Kronecker's delta; and
$u^{S}_{i;L}(x)$ and $u^{S}_{i;T}(x)$ stand for the longitudinal 
and transverse wave components, respectively, of $u^{S}_{i}(x)$.
Also, $k_L^{(m)} = \omega/c_{L}^{(m)}$ is the wavenumber of a longitudinal wave in $\Omega_{m}$,
where $c_L^{(m)} = \sqrt{(\lambda^{(m)} + 2\mu^{(m)})/\rho^{(m)}}$ is the velocity of a longitudinal wave  in $\Omega_{m}$.
Similarly, $k_T^{(m)}$ is the wavenumber of a transverse wave in $\Omega_{m}$,
represented as $k_T^{(m)} = \omega/c_T^{(m)}$, where $c_T^{(m)} = \sqrt{\mu^{(m)}/\rho^{(m)}}$ is the velocity of a transverse wave in $\Omega_{m}$.

\subsection{Boundary integral equations} \label{sec:bie}
Before describing the boundary integral equations, we formulate some layer potential operators.
First, we define the single-layer potential operators $U_{ij}^{(m)}$ with respect to $\Omega_{m}$ as
\begin{align}
  U_{ij}^{(m)} v (x) := \int_{\Gamma} G_{ij}^{(m)} (x, y) v(y) \dd s(y), \quad i, j = 1, 2, \quad x \in \Gamma,
  \label{eq:slp}
\end{align}
where $G_{ij}$ is the fundamental solution of elastodynamics in two dimensions defined by
\begin{multline}
  G_{ij}^{(m)} (x, y) = \frac{\mathrm{i}}{4 \mu^{(m)}} \biggl\{ H_{0}^{(1)}(k_{T}^{(m)} |x - y|) \delta_{ij} \\
  + \frac{1}{(k_{T}^{(m)})^{2}} \left( H_{0}^{(1)}(k_{T}^{(m)} |x - y|) - H_{0}^{(1)}(k_{L}^{(m)} |x - y|) \right)_{,ij} \biggr\}, \quad x, y \in \mathbb{R}^2, \label{eq:funda}
\end{multline}
for $i, j = 1, 2$.
In \eqref{eq:funda}, $H_{0}^{(1)}$ is the zero-th order Hankel function of the first kind.
Then, we define the traction operator $N_{ij}^{(m)} [ \nu(x) ]$ with respect to $\nu(x)$ as
\begin{equation}
  N_{ij}^{(m)} [ \nu(x) ] := \sum_{p, q = 1}^{2} C_{ipjq}^{(m)} \nu_{p} (x) \frac{\partial}{\partial x_{q}}, \quad i, j = 1, 2.
\end{equation}
The traction operator is used to define the double-layer potential operators $T_{ij}^{(m)}$ as
\begin{align}
  T_{ij}^{(m)} v (x) := \mathrm{v.p.} \int_{\Gamma} \qty( \sum_{k = 1}^{2} N_{jk}^{(m)}[\nu(y)] G_{ik}^{(m)} (x, y)) v(y) \dd s(y), \quad i, j = 1, 2, \label{eq:dlp}
\end{align}
where the integral with v.p.\@ is the principal value integral.
The traction of the single-layer potential operators $T_{ij}^{* (m)}$ is defined as
\begin{align}
  T_{ij}^{* (m)} v (x) := \mathrm{v.p.} \int_{\Gamma} \qty( \sum_{k = 1}^{2} N_{ik}^{(m)}[\nu(x)] G_{kj}^{(m)} (x, y)) v(y) \dd s(y), \quad i, j = 1, 2. \label{eq:d_slp}
\end{align}
Finally, the traction of the double-layer potential operators $W_{ij}^{(m)}$ is defined as
\begin{multline}
  W_{ij}^{(m)} v (x) := \mathrm{p.f.} \int_{\Gamma} \qty{\sum_{k = 1}^{2} N_{ik}^{(m)}[\nu(x)] \qty(\sum_{p = 1}^{2} N_{jp}^{(m)}[\nu(y)] G_{kp} (x, y))} v(y) ds(y), \\   i, j = 1, 2,
  \label{eq:d_dlp}
\end{multline}
where the integral with p.f.\@ is the finite part of the divergent integral.
Let $\bm{U}^{(m)}$ be a gathering of the single-layer potential operators as
\begin{equation}
  \bm{U}^{(m)} 
  =
  \mqty(
  U_{11}^{(m)} & U_{12}^{(m)} \\
  U_{21}^{(m)} & U_{22}^{(m)}
  ). \label{eq:gatherd_lp}
\end{equation}
We define $\bm{T}^{(m)}$, $\bm{T}^{* (m)}$, and $\bm{W}^{(m)}$ similarly.
Then, a combination of the single-layer potential operators and appropriate functions $v_1(x)$ and $v_2(x)$ expressed as
\begin{equation}
  \bm{U}^{(m)} 
  \mqty(
  v_1 (x) \\
  v_2 (x)
  )
  =
  \mqty(
  U_{11}^{(m)} & U_{12}^{(m)} \\
  U_{21}^{(m)} & U_{22}^{(m)}
  )
  \mqty(
  v_1 (x) \\
  v_2 (x)
  )
  \label{eq:set_of_lp}
\end{equation}
is called the single-layer potential (with respect to $\Omega_m$).
For the other three block operators,
we similarly refer to the same relations as the double-layer potential,
the traction of single-layer potential, and the traction of the double-layer potential, respectively.

We now describe the boundary integral equations corresponding
to the transmission problems \eqref{eq:navie}--\eqref{eq:raditation} for elastic waves.
The PMCHWT formulation \cite{POGGIO1973159, Chang1977surface} is one of the most famous boundary integral equation formulations
for transmission problems.
In the field of boundary element methods, there is a well-known issue that the existence of
frequencies called fictitious eigenvalues \cite{CHEN1998529, misawa2017boundary}
cause a loss of the uniqueness of the solution even when the original boundary value problem is uniquely solvable.
It is important to employ boundary integral equations that prevent these fictitious eigenvalues from appearing at real frequencies.
The PMCHWT formulation is a formulation of the boundary integral equations expressed as
  \begin{multline}
    \mqty(
    -\qty( U_{11}^{(0)} +       U_{11}^{(1)}) & -\qty( U_{12}^{(0)} + U_{12}^{(1)}) &      T_{11}^{(0)} + T_{11}^{(1)} & T_{12}^{(0)} + T_{12}^{(1)} \\
    -\qty( U_{21}^{(0)} +       U_{21}^{(1)}) & -\qty( U_{22}^{(0)} + U_{22}^{(1)}) &      T_{21}^{(0)} + T_{21}^{(1)} & T_{22}^{(0)} + T_{22}^{(1)} \\
    -\qty(T_{11}^{*(0)} + T_{11}^{*(1)}) & -\qty(T_{12}^{*(0)} + T_{12}^{*(1)}) & W_{11}^{(0)} + W_{11}^{(1)} & W_{12}^{(0)} + W_{12}^{(1)}\\
    -\qty(T_{21}^{*(0)} + T_{21}^{*(1)}) & -\qty(T_{22}^{*(0)} + T_{22}^{*(1)}) & W_{21}^{(0)} + W_{21}^{(1)} & W_{22}^{(0)} + W_{22}^{(1)}
    )
    \mqty(
    t_{1}(x) \\
    t_{2}(x) \\
    u_{1}(x) \\
    u_{2}(x) 
    )
    \\
    =
    \mqty(
         {u_{1}^{I}}(x) \\
         {u_{2}^{I}}(x) \\
         {t_{1}^{I}}(x) \\
         {t_{2}^{I}}(x) 
         ),
         \quad x \in \Gamma,
         \label{eq:pmchwt}
  \end{multline}
where $t_{i}^{I} (x) := \sum_{j = 1}^{2} N_{ij}^{(m)}[\nu(x)] u_{j}^{I}(x)$,
which has no fictitious eigenvalue on real frequencies.
By using \eqref{eq:gatherd_lp}, they can be simplified to
\begin{equation}
  \mqty(
  -\qty(\bm{U}^{(0)} + \bm{U}^{(1)})   & \bm{T}^{(0)} + \bm{T}^{(1)} \\
  -\qty(\bm{T}^{*(0)} + \bm{T}^{* (1)}) & \bm{W}^{(0)} + \bm{W}^{(1)}
  )
  \mqty(
  t(x) \\
  u(x)
  )
  =
  \mqty(
       {u^{I}}(x) \\
       {t^{I}}(x) \\
       ),
       \quad x \in \Gamma, \label{eq:blocked_pmchwt_bie}
\end{equation}
where $t (x) := (t_1 (x), t_2 (x))$ and $t^I (x) := (t_1^I (x), t_2^I (x))$.
We see that the PMCHWT formulated boundary integral equation has eight (block) layer potentials from \eqref{eq:blocked_pmchwt_bie}.

Another formulation that can avoid the fictitious eigenvalue problem is the boundary integral equations
employing the Burton--Miller method \cite{burton1971application} expressed as
  \begin{multline}
    \mqty(
    T_{11}^{(0)} + \frac{1}{2} + \alpha W_{11}^{(0)} & T_{12}^{(0)} + \alpha W_{12}^{(0)} & -U_{11}^{(0)} - \alpha \qty(T_{11}^{* (0)} - \frac{1}{2}) & -U_{12}^{(0)} - \alpha T_{12}^{* (0)}      \\
    T_{21}^{(0)} + \alpha W_{21}^{(0)} & T_{22}^{(0)} + \frac{1}{2} + \alpha W_{22}^{(0)} & -U_{21}^{(0)} - \alpha T_{21}^{* (0)} & -U_{22}^{(0)} - \alpha \qty( T_{22}^{* (0)} - \frac{1}{2})   \\
    T_{11}^{(1)} - \frac{1}{2} & T_{12}^{(1)} & -U_{11}^{(1)} & -U_{12}^{(1)} \\
    T_{21}^{(1)} & T_{22}^{(1)} - \frac{1}{2} & -U_{21}^{(1)} & -U_{22}^{(1)}
    ) \\
    \times
    \mqty(
    u_{1}(x) \\
    u_{2}(x) \\
    t_{1}(x) \\
    t_{2}(x)
    )
    =
    \mqty(
         {u_{1}^{I}}(x) + \alpha {t_{1}^{I}}(x) \\
         {u_{2}^{I}}(x) + \alpha {t_{2}^{I}}(x) \\
         0\\
         0
         ),
         \quad x \in \Gamma.
         \label{eq:bm}
  \end{multline}
Here, $\alpha$ is a constant in the Burton--Miller method satisfying $\Im (\alpha) \neq 0$.
In this study, $\alpha = \mathrm{i}/(\mu^{(0)}k_T^{(0)})$ is used.
We call \eqref{eq:bm} the Burton--Miller formulation in this work.
Similarly to the PMCHWT formulation, it can be simplified to
\begin{equation}
  \mqty(
  \qty(\bm{T}^{(0)} + \frac{1}{2} I_2 + \alpha \bm{W}^{(0)})   & -\bm{U}^{(0)} - \alpha \qty(\bm{T}^{*(0)} - \frac{1}{2} I_2) \\
  \bm{T}^{(1)} - \frac{1}{2} I_2 & - \bm{U}^{(1)}
  )
  \mqty(
  u(x) \\
  t(x)
  )
  =
  \mqty(
  u^{I}(x) + \alpha t^{I}(x) \\
  0
  ),
  \quad x \in \Gamma, \label{eq:blocked_bm_bie}
\end{equation}
where $I_2$ is an identity matrix of size $2 \times 2$.
A numerical method based on the Burton--Miller formulation \eqref{eq:bm} is expected to be faster than the PMCHWT formulation \eqref{eq:blocked_pmchwt_bie} because the number of (block) layer potentials of the Burton--Miller formulation is six while that in the PMCHWT is eight.

\begin{remark}
Formulations \eqref{eq:pmchwt} and \eqref{eq:bm} are suitable for fast direct solvers
because they deal with the (physical) direct unknowns.
When we use some indirect formulations, we also need fast matrix-vector multiplication methods to obtain the (physical) unknowns on the boundary.
Such fast matrix-vector multiplication methods can be constructed \cite{YESYPENKO2025113707, matsumoto2024fast} within the framework of the proxy method discussed in this work,
but it is inefficient because it requires unnecessary calculations when using physical unknowns.
\end{remark}

\begin{remark}
In this study, we discuss the fast method up to solving the boundary integral equation \eqref{eq:pmchwt} or \eqref{eq:bm}.
After the solutions $u_{i}$ and $t_{i}$ on the boundary $\Gamma$ are obtained,
the solution of the boundary value problems \eqref{eq:navie}--\eqref{eq:raditation} in domain $\Omega_{0}$ can be obtained by using the following integral representation:
\begin{multline}
  u_i(x) = u^{I}_{i}(x) - \int_{\Gamma} \qty(N_{jk}^{(0)}[\nu(y)] G_{ik}^{(0)} (x, y)) u_{j}(y) \dd s(y) + \int_{\Gamma} G_{ij}^{(0)} (x, y) t_{j}(y) \dd s(y), \\ x \in \Omega_{0}.
\end{multline}
Similarly, the solution of the boundary value problems \eqref{eq:navie}--\eqref{eq:raditation} in domain $\Omega_{1}$ can be obtained by using the following integral representation:
\begin{align}
  u_i(x) = \int_{\Gamma} \qty(N_{jk}^{(1)}[\nu(y)] G_{ik}^{(1)} (x, y)) u_{j}(y) \dd s(y) - \int_{\Gamma} G_{ij}^{(1)} (x, y) t_{j}(y) \dd s(y), \quad x \in \Omega_{1}.
\end{align}
\end{remark}

\section{Discretization} \label{sec:Discretization}
As mentioned in Section \ref{sec:intro}, the versatility
of the proposed method is emphasized in this work.
A discretization technique that fits this idea is discussed.

\subsection{Galerkin method} \label{sec:galerkin}

At first, the boundary $\Gamma$ is approximated as a $N_n$ polygon,
where $N_n$ is the number of nodes.
Let $\{ x^{s} \}_{s = 1}^{N_n}$ be
a set of nodes whose coordinates are the vertices of the $N_n$ polygon.
Let $N_e$ be the number of elements of the $N_n$ polygon.
In this work, $N_e = N_n =: N$ because $\Gamma \subset \mathbb{R}^2$
is a closed Jordan curve and straight line segments are used for boundary elements.
The case $N_e \neq N_n$ arises from the scattering analysis by cracks (open arcs) even in two dimensions.
Let $\{ \varphi^{s} (x) \}_{s = 1}^{N_n}$ and $\{ \psi^{s} (x) \}_{s = 1}^{N_e}$ be
a set of the piecewise linear bases
and a set of the piecewise constant bases
on the approximated $\Gamma$, respectively.
By using $\varphi^{s}$, the unknown functions $u_i$ of \eqref{eq:pmchwt} or \eqref{eq:bm}
on the boundary are approximated as
\begin{align}
  u_{i} (x) \approx \sum_{s = 1}^{N_n} \varphi^{s} (x) u_{i}^{s},
\end{align}
where $u^{s} = (u_{1}^{s}, u_{2}^{s}) = (u_{1} (x^s), u_{2} (x^s)) \in \mathbb{C}^2$.
Similarly, the unknown functions $t_i$ of \eqref{eq:pmchwt} or \eqref{eq:bm}
on the boundary are approximated by using $\psi^{s}$ as
\begin{align}
  t_{i} (x) \approx \sum_{s = 1}^{N_e} \psi^{s} (x) t_{i}^{s},
\end{align}
where $t^{s} = (t_{1}^{s}, t_{2}^{s}) = (t_{1} (x_c^s), t_{2} (x_c^s)) \in \mathbb{C}^2$.
Here, $x_c^s$ is the centroid of $s$-th element.
Note that the displacement $u_i$ is approximated by piecewise linear bases
while the traction $t_i$ is approximated by piecewise constant bases,
because they have different smoothness on $\Gamma$.

When using these basis assemblages, an appropriate selection of test functions for the PMCHWT formulation \eqref{eq:pmchwt} can be outlined as shown in Table \ref{tb:pmchwt}.
Similarly, an appropriate selection of test functions for the Burton--Miller formulation \eqref{eq:bm} can be outlined as shown in Table \ref{tb:bm}.
By selecting the test functions as described above,
we can make the coefficient matrices for the discrete systems for both formulations square even in the general case.
\begin{table}[h]
  \centering
  \caption{Assembly of test and basis functions for the PMCHWT formulation. $\varphi$ and $\psi$ are used as generic notations for the set of piecewise linear bases $\{ \varphi^{s} (x) \}_{s = 1}^{N_n}$ and the set of piecewise constant bases $\{ \psi^{s} (x) \}_{s = 1}^{N_e}$, respectively}
  \label{tb:pmchwt}
  \begin{tabular}{c|cc}
    Test function & $\, \psi$ is used for approximation of $t$ & $\quad \varphi$ is used for approximation of $u$ \\
    \hline
    $ \psi$ & $-(\bm{U}^{(0)} + \bm{U}^{(1)})$             & $\bm{T}^{(0)} + \bm{T}^{(1)}$ \\
    $\varphi$    & $-({\bm{T}^{\ast (0)}} + {\bm{T}^{\ast (1)}})$ & $\bm{W}^{(0)} + \bm{W}^{(1)}$
  \end{tabular}
\end{table}
\begin{table}[h]
  \centering
  \caption{Assembly of test and basis functions for the Burton--Miller formulation. $\varphi$ and $\psi$ are used as generic notations for the set of piecewise linear bases $\{ \varphi^{s} (x) \}_{s = 1}^{N_n}$ and the set of piecewise constant bases $\{ \psi^{s} (x) \}_{s = 1}^{N_e}$, respectively}
  \label{tb:bm}
  \begin{tabular}{c|cc}
   Test function & $\, \varphi$ is used for approximation of $u$ & $\quad \psi$ is used for approximation of $t$ \\
    \hline
    $\varphi$ & $\bm{T}^{(0)} + (1/2) I_2 + \alpha \bm{W}^{(0)}$             & $- \bm{U}^{(0)} - \alpha \qty( \bm{T}^{* (0)} - (1/2) I_2)$ \\
    $\psi$    & $\bm{T}^{(1)} - (1/2) I_2$ & $-\bm{U}^{(1)}$
  \end{tabular}
\end{table}

We now precisely describe how these functions are used. By using basis functions and the $r$-th test functions $\psi^{r}$ and $\varphi^{r}$, the
PMCHWT formulation is discretized as the following system of equations:
\begin{multline}
  \int_\Gamma \psi^{r}(x) \qty{ \qty( -\sum_{j =1}^{2}\qty(U_{ij}^{(0)} + U_{ij}^{(1)}) \sum_{s=1}^{N_e}\psi^{s}(x) t_{j}^{s})
    + \qty( \sum_{j = 1}^{2} \qty(T_{ij}^{(0)} + T_{ij}^{(1)}) \sum_{s=1}^{N_n} \varphi^{s}(x) u_{j}^{s})
  } \dd s(x) \\
  = \int_{\Gamma} \psi^{r}(x) u^{I}_{i}(x) \dd s(x), \quad i = 1, 2, \quad r = 1, 2, \ldots, N_{e}, \quad x \in \Gamma,
  \label{eq:gal_pm1}
\end{multline}
\begin{multline}
  \int_\Gamma \varphi^{r}(x) \qty{ \qty( -\sum_{j =1}^{2}\qty(T_{ij}^{*(0)} + T_{ij}^{*(1)}) \sum_{s=1}^{N_e}\psi^{s}(x) t_{j}^{s})
    + \qty( \sum_{j = 1}^{2} \qty(W_{ij}^{(0)} + W_{ij}^{(1)}) \sum_{s=1}^{N_n} \varphi^{s}(x) u_{j}^{s})
  } \dd s(x) \\
  = \int_{\Gamma} \varphi^{r}(x) t^{I}_{i}(x) \dd s(x), \quad i = 1, 2, \quad r = 1, 2, \ldots, N_{n}, \quad x \in \Gamma.
  \label{eq:gal_pm2}
\end{multline}
Similarly, the Burton--Miller formulation is discretized by the basis functions and the $r$-th test functions $\psi^{r}$ and $\varphi^{r}$ as the following system of equations:
\begin{multline}
  \int_\Gamma \varphi^{r}(x) \left\{ \qty( \sum_{j =1}^{2} \qty(T_{ij}^{(0)} + \frac{\delta_{ij}}{2} + \alpha W_{ij}^{(0)}) \sum_{s=1}^{N_n}\varphi^{s}(x) u_{j}^{s} ) \right. \\
    \left. + \qty( \sum_{j = 1}^{2} \qty{-U_{ij}^{(0)} - \alpha \qty(T_{ij}^{*(0)} - \frac{\delta_{ij}}{2})} \sum_{s=1}^{N_e} \psi^{s}(x) t_{j}^{s})
    \right\} \dd s(x) \\
  = \int_{\Gamma} \varphi^{r}(x)\qty(u_{i}^{I}(x) + \alpha t^{I}_{i}(x)) \dd s(x), \quad i = 1, 2, \quad r = 1, 2, \ldots, N_{n}, \quad x \in \Gamma.
  \label{eq:gal_bm1}
\end{multline}
\begin{multline}
  \int_\Gamma \psi^{r}(x) \qty{ \qty( \sum_{j =1}^{2}\qty(T_{ij}^{(1)} - \frac{\delta_{ij}}{2}) \sum_{s=1}^{N_n}\varphi^{s}(x) u_{j}^{s})
    - \qty( \sum_{j = 1}^{2} U_{ij}^{(1)} \sum_{s=1}^{N_e} \psi^{s}(x) t_{j}^{s})
  } \dd s(x) \\
  = 0, \quad i = 1, 2, \quad r = 1, 2, \ldots, N_{e}, \quad x \in \Gamma.
  \label{eq:gal_bm2}
\end{multline}
The simplest way to find the discretized solutions $u$ and $t$ on $\Gamma$ in this Galerkin discretization is to solve the system of \eqref{eq:gal_pm1} and \eqref{eq:gal_pm2} for the PMCHWT formulation and the system of \eqref{eq:gal_bm1} and \eqref{eq:gal_bm2} for the Burton--Miller formulation,
which are together treated as the conventional boundary element method.

\begin{remark}
We employ a regularization technique from Appendix L.5 of \cite{YosK:2001}
to compute the discretized $W_{ij}^{(m)}$, similarly to \cite{MATSUMOTO2025106148}.
The regularization technique used in \cite{YosK:2001}
transfers the derivatives of the fundamental solution to the basis functions
using integration by parts, which can be easily modified for two-dimensional cases.
\end{remark}

\subsection{Block formulation}
The proposed fast direct solver is mainly based on a low-rank approximation of off-diagonal parts of a matrix.
We transform the (discretized) linear equations into block forms in preparation for a low-rank approximation.
The algorithm for the proposed fast direct solver in this work is almost the same as that in \cite{MATSUMOTO2025106148},
except that the tree partitioning is performed separately
for the piecewise linear bases $\{ \varphi^s \}$ and the piecewise constant bases $\{ \psi^s \}$.
In the implementation, cells of the tree for nodes $\{ J_i^{\varphi} \}$ and for elements $\{ J_i^\psi \}$ are constructed,
because $\{ \varphi^s \}$ and $\{ \psi^s \}$ correspond to the coordinates of nodes and element centers, respectively.
We assume that the structures of the tree for nodes and elements are the same in this work.

Let $\{ \varphi^s \}_{s = 1}^{N_n}$ and $\{ \psi^s \}_{s = 1}^{N_e}$
be partitioned to cells $\{J_1^{\varphi}, J_2^{\varphi}, J_3^{\varphi}, \ldots, J_p^{\varphi} \}$
and $\{J_1^\psi, J_2^\psi, J_3^\psi, \ldots, J_p^\psi \}$
according to a binary tree at the leaf level, where $p$ is the number of cells at the leaf level.
The cells $J_{i}^{\varphi}$ and $J_{i}^{\psi}$ are sequences of indices
such that the $s$-th indices in $J_{i}^{\varphi}$ and $J_{i}^{\psi}$ are expressed as
\begin{equation}
  J_{i}^{\varphi}(s) = s + \qty(\sum_{k = 1}^{i - 1} |J_{k}^{\varphi}|), \quad J_{i}^{\psi}(s) = s + \qty(\sum_{k = 1}^{i - 1} |J_{k}^{\psi}|),
\end{equation}
for $i = 1, 2, \ldots, p$, where $| J_{i}^{\varphi} |$ and $| J_{i}^{\psi} |$ are the number of elements in the sequences $ J_{i}^{\varphi} $ and $ J_{i}^{\psi} $, respectively, 
and the summation for $i = 1$ in the above relations is the empty sum.
The system of the discretized PMCHWT or Burton--Miller formulation,
with coefficient matrix $A$, solution vector $x$, and right-hand side vector $f$,
is partitioned into blocks by cells of the tree.
Here, the size of $A$, $x$, and $f$ are $4N \times 4N$, $4N$, and $4N$, respectively.
In addition, to perform a low-rank approximation in the vector field appropriately,
the matrix components corresponding to the interactions between each cell are 
grouped by the directional components (in the vector field) along which the layer potential operators act, as shown in \eqref{eq:zblock}.
Then, the system is written in the following block form:
\begin{equation}
  \mqty(
  A_{11} & A_{12} & \cdots & A_{1p} \\
  A_{21} & A_{22} & \cdots & A_{2p} \\
  \vdots & \vdots &\ddots &\vdots \\
  A_{p1} & A_{p2} & \cdots & A_{pp}
  )
  \mqty(
  x_1 \\
  x_2 \\
  \vdots \\
  x_p
  )
  =
  \mqty(
  f_1 \\
  f_2 \\
  \vdots \\
  f_p
  ),
  \label{eq:block_linear}
\end{equation}
where each $A_{ij}$ corresponds to the interactions between cells
$J_{i}^{\varphi}$, $J_{i}^{\psi}$, $J_{j}^{\varphi}$, and $J_{j}^{\psi}$ for $i, j = 1, 2, \ldots, p$.
Similarly, $x_{i}$ corresponds to the solution on the nodes and elements in $J_{i}^{\varphi}$ and $J_{i}^{\psi}$,
and $f_{i}$ corresponds to the discretized incident wave in $J_{i}^{\varphi}$ and $J_{i}^{\psi}$.
More precisely, $A_{ij}$ in the PMCHWT formulation is defined by
\begin{equation}
  A_{ij} =
    \mqty(
    -\qty( \qty[\bm{U}^{(0)}]^{J_{i}^{\psi} J_{j}^{\psi}} +       \qty[\bm{U}^{(1)}]^{J_{i}^{\psi} J_{j}^{\psi}}) & \qty[\bm{T}_{}^{(0)}]^{J_{i}^{\psi} J_{j}^{\varphi}} + \qty[\bm{T}_{}^{(1)}]^{J_{i}^{\psi} J_{j}^{\varphi}} \\
    -\qty(\qty[\bm{T}_{}^{*(0)}]^{J_{i}^{\varphi} J_{j}^{\psi}} + \qty[\bm{T}_{}^{*(1)}]^{J_{i}^{\varphi} J_{j}^{\psi}}) & \qty[\bm{W}_{}^{(0)}]^{J_{i}^{\varphi} J_{j}^{\varphi}} + \qty[\bm{W}_{}^{(1)}]^{J_{i}^{\varphi} J_{j}^{\varphi}} \\
    ),
  \label{eq:pmchwt_discrete}
\end{equation}
while it in the Burton--Miller formulation is defined by
\begin{equation}
  A_{ij} =
  \mqty(
  \qty[\bm{T}^{(0)}]^{J_{i}^{\varphi} J_{j}^{\varphi}} + \frac{I^{J_{i}^{\varphi} J_{j}^{\varphi}}}{2} + \alpha \qty[\bm{W}^{(0)}]^{J_{i}^{\varphi} J_{j}^{\varphi}}   & \quad -\qty[\bm{U}^{(0)}]^{J_{i}^{\varphi} J_{j}^{\psi}} - \alpha \qty(\qty[\bm{T}^{*(0)}]^{J_{i}^{\varphi} J_{j}^{\psi}} - \frac{I^{J_{i}^{\varphi} J_{j}^{\psi}}}{2}) \\
    \qty[\bm{T}^{(1)}]^{J_{i}^{\psi} J_{j}^{\varphi}} - \frac{I^{J_{i}^{\psi} J_{j}^{\varphi}}}{2}  & - \qty[\bm{U}^{(1)}]^{J_{i}^{\psi} J_{j}^{\psi}}
  ).
  \label{eq:bm_discrete}
\end{equation}
Here, the notation $\qty[\bm{Z}^{(m)}]^{J_{i}^{\zeta} J_{j}^{\xi}}$ indicates
a matrix of size $2|J_{i}^{\zeta}| \times 2|J_{j}^{\xi}|$ defined by
\begin{equation}
  \qty[\bm{Z}^{(m)}]^{J_{i}^{\zeta} J_{j}^{\xi}} :=
  \mqty(
  \qty[Z_{11}^{(m)}]^{J_{i}^{\zeta} J_{j}^{\xi}} & \qty[Z_{12}^{(m)}]^{J_{i}^{\zeta} J_{j}^{\xi}} \\
  \qty[Z_{21}^{(m)}]^{J_{i}^{\zeta} J_{j}^{\xi}} & \qty[Z_{22}^{(m)}]^{J_{i}^{\zeta} J_{j}^{\xi}}
  ), \label{eq:zblock}
\end{equation}
where $Z_{ij}^{(m)}$ is an alternative for $U_{ij}^{(m)}$, $T_{ij}^{(m)}$, $T_{ij}^{*(m)}$, and $W_{ij}^{(m)}$.
Also, $\zeta$ or $\xi$ is an alternative for $\varphi$ and $\psi$.
Let $\qty[Z_{pq}^{(m)}]_{rs}^{J_{i}^{\zeta} J_{j}^{\xi}}$ be the $(r, s)$ component of $\qty[Z_{pq}^{(m)}]^{J_{i}^{\zeta} J_{j}^{\xi}}$ defined as
\begin{equation}
  \qty[Z_{pq}^{(m)}]_{rs}^{J_{i}^{\zeta} J_{j}^{\xi}} :=
  \int_{\Gamma} \zeta^{J_{i}^{\zeta} (r)}(x) Z_{pq}^{(m)}\xi^{J_{j}^{\xi}(s)}(x) \dd s(x).
\end{equation}
In \eqref{eq:bm_discrete}, $I^{J_{i}^{\zeta} J_{j}^{\xi}}$ stands for
the matrix with respect to a set of test functions $\{ \zeta^s \}_{s=1}^{|J_{i}^{\zeta}|}$
and a set of basis functions $\{ \xi^{s} \}_{s=1}^{|J_{j}^{\xi}|}$ defined by
\begin{equation}
  I^{J_{i}^{\zeta} J_{j}^{\xi}} :=
  \mqty(
  I_{0}^{J_{i}^{\zeta} J_{j}^{\xi}} & 0 \\
  0 & I_{0}^{J_{i}^{\zeta} J_{j}^{\xi}}
  ),
\end{equation}
where the $(r, s)$ component of $I_{0}^{J_{i}^{\zeta} J_{j}^{\xi}}$ is defined by
\begin{equation}
  \qty[ I_{0}^{J_{i}^{\zeta} J_{j}^{\xi}}]_{rs} = \int_\Gamma \zeta^{J_{i}^{\zeta}(r)}(x) \xi^{J_{j}^{\xi}(s)}(x) \dd s(x).
\end{equation}
Again, $\zeta$ or $\xi$ is an alternative for $\varphi$ and $\psi$.
For example, if we set $\zeta = \varphi$ and $\xi = \psi$, the resulting relation becomes
\begin{equation}
  \qty[I_{0}^{J_{i}^{\varphi} J_{j}^{\psi}}]_{rs} = \int_\Gamma \varphi^{J_{i}^{\varphi}(r)}(x) \psi^{J_{j}^{\psi}(s)}(x) \dd s(x).
\end{equation}
In \eqref{eq:block_linear}, $x_i \in \mathbb{C}^{2|J_{i}^{\varphi}| + 2|J_{i}^{\psi}|}$ in the PMCHWT formulation is defined as
\begin{equation}
  x_i =
    \mqty(
    x_{t_1}^{J_{i}^{\psi}} \\
    x_{t_2}^{J_{i}^{\psi}} \\
    x_{u_1}^{J_{i}^{\varphi}} \\
    x_{u_2}^{J_{i}^{\varphi}}
    ),
\end{equation}
while it in the Burton--Miller formulation is defined as
\begin{equation}
  x_i =
    \mqty(
    x_{u_1}^{J_{i}^{\varphi}} \\
    x_{u_2}^{J_{i}^{\varphi}} \\
    x_{t_1}^{J_{i}^{\psi}} \\
    x_{t_2}^{J_{i}^{\psi}} \\
    ),
\end{equation}
where the vectors $x_{u_j}^{J_{i}^{\varphi}} \in \mathbb{C}^{|J_{i}^{\varphi}|}$ and 
$x_{t_j}^{J_{i}^{\psi}} \in \mathbb{C}^{|J_{i}^{\psi}|}$ have an $r$-th component expressed as
\begin{equation}
  \qty[x_{u_j}^{J_{i}^{\varphi}}]_r := u_{j}^{J_{i}^{\varphi}(r)}, \quad j = 1, 2, \quad r = 1, 2, \ldots, |J_{i}^{\varphi}|,
\end{equation}
and
\begin{equation}
  \qty[x_{t_j}^{J_{i}^{\psi}}]_r := t_{j}^{J_{i}^{\psi}(r)}, \quad j = 1, 2, \quad r = 1, 2, \ldots, |J_{i}^{\psi}|.
\end{equation}
Similarly, $f_i \in \mathbb{C}^{2|J_{i}^{\varphi}| + 2|J_{i}^{\psi}|}$ in the PMCHWT formulation is defined as
\begin{equation}
  f_i =
    \mqty(
    f_{u_1}^{J_{i}^{\psi}} \\
    f_{u_2}^{J_{i}^{\psi}} \\
    f_{t_1}^{J_{i}^{\varphi}} \\
    f_{t_2}^{J_{i}^{\varphi}}
    ),
\end{equation}
while it in the Burton--Miller formulation is defined as
\begin{equation}
  f_i =
    \mqty(
    f_{u_1}^{J_{i}^{\varphi}} + \alpha f_{t_1}^{J_{i}^{\varphi}} \\
    f_{u_2}^{J_{i}^{\varphi}} + \alpha f_{t_2}^{J_{i}^{\varphi}} \\
    0 \\
    0
    ),
    \label{eq:rhs_f_of_BM}
\end{equation}
where the vectors $f_{u_j}^{J_{i}^{\varphi}} \in \mathbb{C}^{|J_{i}^{\varphi}|}$,
$f_{t_j}^{J_{i}^{\varphi}} \in \mathbb{C}^{|J_{i}^{\varphi}|}$,
and $f_{u_j}^{J_{i}^{\psi}} \in \mathbb{C}^{|J_{i}^{\psi}|}$
have their $r$-th components expressed as
\begin{align}
  \qty[f_{u_j}^{J_{i}^{\varphi}}]_r &= \int_{\Gamma} \varphi^{J_i^\varphi (r)}(x) u^{I}_{j}(x) \dd s(x), \quad j = 1, 2, \quad r = 1, 2, \ldots, |J_{i}^{\varphi}|, \\
  \qty[f_{t_j}^{J_{i}^{\varphi}}]_r &= \int_{\Gamma} \varphi^{J_i^\varphi (r)}(x) t^{I}_{j}(x) \dd s(x), \quad j = 1, 2, \quad r = 1, 2, \ldots, |J_{i}^{\varphi}|, \\
  \qty[f_{u_j}^{J_{i}^{\psi}}]_r &= \int_{\Gamma} \psi^{J_i^\psi (r)}(x) u^{I}_{j}(x) \dd s(x), \quad j = 1, 2, \quad r = 1, 2, \ldots, |J_{i}^{\psi}|,
\end{align}
respectively.
In \eqref{eq:rhs_f_of_BM}, each $0$ means a zero vector of size $|J_{i}^{\psi}|$.

\section{Proposed fast method} \label{sec:proposed}
In this paper, only the single-level method for the proposed fast direct solver is formulated.
The multi-level version can be formulated easily in the same way as done for elastic wave scattering by a cavity \cite{MATSUMOTO2025106148}.
In a multi-level algorithm, note that skeletons $S_{ij}$ in \eqref{eq:a=usv} are only required at the top level
while $L_i$ and $R_i$ in \eqref{eq:a=usv} are required for all levels.

\subsection{Low-rank approximations via the proxy method}
It is known that the off-diagonal block $A_{ij}$ $(i \neq j)$ of \eqref{eq:block_linear}
can be low-rank approximated in formulations of boundary element methods by the property of fundamental solutions.
Here, we discuss the weak admissibility type low-rank approximation for $A_{ij} \in \mathbb{C}^{2 \qty(|J_i^{\varphi}| + |J_i^\psi|) \times 2 \qty(|J_i^{\varphi}| + |J_i^\psi|)}$.
In this section, we assume $n = |J_i^{\varphi}| = |J_i^{\psi}|$ $(i = 1, 2, \ldots, p)$ for simplicity,
so that $A_{ij}$ can be expressed as $A_{ij} \in \mathbb{C}^{4n \times 4n}$.
Then, $A_{ij}$ is low-rank approximated to
\begin{equation}
  A_{ij} = L_{i} S_{ij} R_{j}, \quad i \neq j,
  \label{eq:a=usv}
\end{equation}
where
$S_{ij} \in \mathbb{C}^{4k \times 4k}$ are the so-called skeletons of $A_{ij}$; $k$ is the number of skeletons; and
$L_i \in \mathbb{C}^{4n \times 4k}$ and $R_i \in \mathbb{C}^{4k \times 4n}$
are left and right coefficients of low-rank approximations, respectively, with respect to $S_{ij}$.
We assume $k \ll n$.
In this section, we assume that $k$ is constant.
It is important for an efficient fast direct solver
that $L_i$ and $R_j$ are shared for the $i$-th row block and $j$-th column block, respectively.

For an accurate description of the skeletons $S_{ij}$,
we introduce the subsequences $\tilde{J}_{i}^{\varphi}$ and $\tilde{J}_{i}^{\psi}$ of $J_{i}^{\varphi}$ and $J_{i}^{\psi}$, respectively.
The method for calculating subsequences $\tilde{J}_{i}^{\varphi}$ and $\tilde{J}_{i}^{\psi}$ and 
coefficients $L_i$ and $R_i$ is described later.
Let $\tilde{J}_{i}^{\varphi}$ and $\tilde{J}_{i}^{\psi}$ have $k$ elements
from $J_{i}^{\varphi}$ and $J_{i}^{\psi}$, respectively.
Then $S_{ij}$ is defined by
\begin{equation}
  S_{ij} := M(\tilde{J}_{i}, \tilde{J}_{j}),
\end{equation}
where $\tilde{J}_i$ is the ordered pair $( \tilde{J}_{i}^{\varphi}, \tilde{J}_{i}^{\psi} )$
and a matrix $M (\tilde{J}_i, \tilde{J}_j)$
and its sub-matrices $M^{t} (\tilde{J}_i, \tilde{J}_j)$, $t = 1, 2, 3, 4$, are defined by
\begin{align}
  M(\tilde{J}_{i}, \tilde{J}_{j}) &:=
  \mqty(
  M^{1}(\tilde{J}_{i}, \tilde{J}_{j}) & M^{2}(\tilde{J}_{i}, \tilde{J}_{j}) \\
  M^{3}(\tilde{J}_{i}, \tilde{J}_{j}) & M^{4}(\tilde{J}_{i}, \tilde{J}_{j})
  ) \\
  &:=
  \begin{cases}
    \begin{split}
    \mqty(
    -\qty( \qty[\bm{U}^{(0)}]^{\tilde{J}_{i}^{\psi} \tilde{J}_{j}^{\psi}} +       \qty[\bm{U}^{(1)}]^{\tilde{J}_{i}^{\psi} \tilde{J}_{j}^{\psi}}) & \qty[\bm{T}^{(0)}]^{\tilde{J}_{i}^{\psi} \tilde{J}_{j}^{\varphi}} + \qty[\bm{T}^{(1)}]^{\tilde{J}_{i}^{\psi} \tilde{J}_{j}^{\varphi}} \\
    -\qty(\qty[\bm{T}^{*(0)}]^{\tilde{J}_{i}^{\varphi} \tilde{J}_{j}^{\psi}} + \qty[\bm{T}^{*(1)}]^{\tilde{J}_{i}^{\varphi} \tilde{J}_{j}^{\psi}}) & \qty[\bm{W}^{(0)}]^{\tilde{J}_{i}^{\varphi} \tilde{J}_{j}^{\varphi}} + \qty[\bm{W}^{(1)}]^{\tilde{J}_{i}^{\varphi} \tilde{J}_{j}^{\varphi}} \\
    ), \\
    \text{for the PMCHWT formulation,} \\ \quad
    \end{split}
    \\
    \begin{split}
  \mqty(
  \qty[\bm{T}^{(0)}]^{\tilde{J}_{i}^{\varphi} \tilde{J}_{j}^{\varphi}} + \alpha \qty[\bm{W}^{(0)}]^{\tilde{J}_{i}^{\varphi} \tilde{J}_{j}^{\varphi}}   & -\qty[\bm{U}^{(0)}]^{\tilde{J}_{i}^{\varphi} \tilde{J}_{j}^{\psi}} - \alpha \qty[\bm{T}^{*(0)}]^{\tilde{J}_{i}^{\varphi} \tilde{J}_{j}^{\psi}} \\
    \qty[\bm{T}^{(1)}]^{\tilde{J}_{i}^{\psi} \tilde{J}_{j}^{\varphi}} & - \qty[\bm{U}^{(1)}]^{\tilde{J}_{i}^{\psi} \tilde{J}_{j}^{\psi}}
  ), \\
  \quad \text{for the Burton--Miller formulation,}
    \end{split}
  \end{cases}
\end{align}
for $i \neq j$ and $i, j = 1, 2, \ldots, p$.
Let $M_{pq}^{t} (\tilde{J}_{i}, \tilde{J}_{j})$, $p, q = 1, 2$ and $i \neq j$, be a submatrix of $M^{t}$, $t = 1, 2, 3, 4$, such that
\begin{equation}
  M^{t}(\tilde{J}_{i}, \tilde{J}_{j}) = 
    \mqty(
    M_{11}^{t}(\tilde{J}_{i}, \tilde{J}_{j}) & M_{12}^{t}(\tilde{J}_{i}, \tilde{J}_{j}) \\
    M_{21}^{t}(\tilde{J}_{i}, \tilde{J}_{j}) & M_{22}^{t}(\tilde{J}_{i}, \tilde{J}_{j})
    ),
\end{equation}
where each $M_{pq}^{t}(\tilde{J}_{i}, \tilde{J}_{j})$ is defined by
\begin{align}
  M_{pq}^{1}(\tilde{J}_{i}, \tilde{J}_{j}) &:=
  \begin{cases}
    -\qty( \qty[U_{pq}^{(0)}]^{\tilde{J}_{i}^{\psi} \tilde{J}_{j}^{\psi}} + \qty[U_{pq}^{(1)}]^{\tilde{J}_{i}^{\psi} \tilde{J}_{j}^{\psi}})
    , \quad \text{for the PMCHWT formulation,} \\
    \qty[T_{pq}^{(0)}]^{\tilde{J}_{i}^{\varphi} \tilde{J}_{j}^{\varphi}} + \alpha \qty[W_{pq}^{(0)}]^{\tilde{J}_{i}^{\varphi} \tilde{J}_{j}^{\varphi}}
    , \quad \text{for the Burton--Miller formulation,}
  \end{cases} \\
  M_{pq}^{2}(\tilde{J}_{i}, \tilde{J}_{j}) &:=
  \begin{cases}
    \qty[T_{pq}^{(0)}]^{\tilde{J}_{i}^{\psi} \tilde{J}_{j}^{\varphi}} + \qty[T_{pq}^{(1)}]^{\tilde{J}_{i}^{\psi} \tilde{J}_{j}^{\varphi}}
    , \quad \text{for the PMCHWT formulation,} \\
    -\qty[U_{pq}^{(0)}]^{\tilde{J}_{i}^{\varphi} \tilde{J}_{j}^{\psi}} - \alpha \qty[T_{pq}^{*(0)}]^{\tilde{J}_{i}^{\varphi} \tilde{J}_{j}^{\psi}}
    , \quad \text{for the Burton--Miller formulation,}
  \end{cases} \\
  M_{pq}^{3}(\tilde{J}_{i}, \tilde{J}_{j}) &:=
  \begin{cases}
    -\qty( \qty[T_{pq}^{*(0)}]^{\tilde{J}_{i}^{\varphi} \tilde{J}_{j}^{\psi}} + \qty[T_{pq}^{*(1)}]^{\tilde{J}_{i}^{\varphi} \tilde{J}_{j}^{\psi}})
    , \quad \text{for the PMCHWT formulation,} \\
    \qty[T_{pq}^{(1)}]^{\tilde{J}_{i}^{\psi} \tilde{J}_{j}^{\varphi}}
    , \quad \text{for the Burton--Miller formulation,}
  \end{cases} \\
  M_{pq}^{4}(\tilde{J}_{i}, \tilde{J}_{j}) &:=
  \begin{cases}
    \qty[W_{pq}^{(0)}]^{\tilde{J}_{i}^{\varphi} \tilde{J}_{j}^{\varphi}} + \qty[W_{pq}^{(1)}]^{\tilde{J}_{i}^{\varphi} \tilde{J}_{j}^{\varphi}}
    , \quad \text{for the PMCHWT formulation,} \\
    -\qty[U_{pq}^{(1)}]^{\tilde{J}_{i}^{\psi} \tilde{J}_{j}^{\psi}}
    , \quad \text{for the Burton--Miller formulation,}
  \end{cases}
\end{align}
respectively.
To formulate the proxy method, the upper, lower, left, and right halves of $M^p (\tilde{J}_{i}, \tilde{J}_{j})$ denoted by $M_{\mathrm{uh}}^{p}(\tilde{J}_{i}, \tilde{J}_{j})$, $M_{\mathrm{dh}}^{p}(\tilde{J}_{i}, \tilde{J}_{j})$, $M_{\mathrm{lh}}^{p}(\tilde{J}_{i}, \tilde{J}_{j})$, and $M_{\mathrm{rh}}^{p}(\tilde{J}_{i}, \tilde{J}_{j})$, respectively, are defined as
\begin{align}
  M_{\mathrm{uh}}^{p}(\tilde{J}_{i}, \tilde{J}_{j}) &:= \mqty(M_{11}^{p}(\tilde{J}_{i}, \tilde{J}_{j}) & M_{12}^{p}(\tilde{J}_{i}, \tilde{J}_{j})) \\
  M_{\mathrm{dh}}^{p}(\tilde{J}_{i}, \tilde{J}_{j}) &:= \mqty(M_{21}^{p}(\tilde{J}_{i}, \tilde{J}_{j}) & M_{22}^{p}(\tilde{J}_{i}, \tilde{J}_{j})) \\
  M_{\mathrm{lh}}^{p}(\tilde{J}_{i}, \tilde{J}_{j}) &:= \mqty(M_{11}^{p}(\tilde{J}_{i}, \tilde{J}_{j}) \\ M_{21}^{p}(\tilde{J}_{i}, \tilde{J}_{j})) \\
  M_{\mathrm{rh}}^{p}(\tilde{J}_{i}, \tilde{J}_{j}) &:= \mqty(M_{12}^{p}(\tilde{J}_{i}, \tilde{J}_{j}) \\ M_{22}^{p}(\tilde{J}_{i}, \tilde{J}_{j})).
\end{align}
Similarly, the pair $J_i = (J_i^{\varphi}, J_i^{\psi})$ is defined in the same way as above. 
\begin{figure}[tb]
  \centering
  \includegraphics[width=0.70\linewidth]{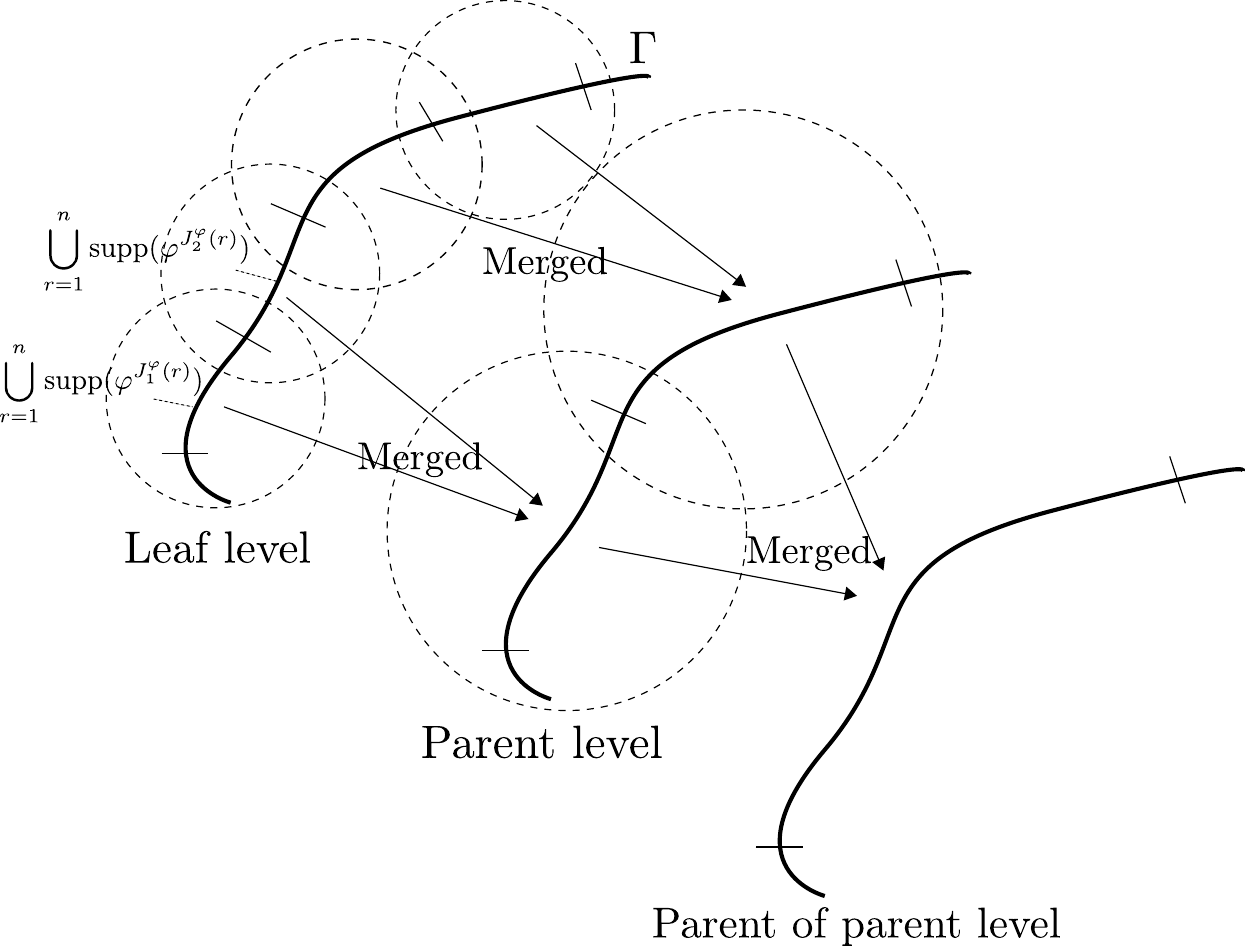}
  \caption{Overview of the proxy method. This figure corresponds to the case of $\{ \varphi^r \}$.
    The circles with dashed lines are proxy boundaries. The boundary $\Gamma$ is partitioned to cells,
    which correspond to the sum of each support of the basis functions
    denoted as $\bigcup_{r=1}^{n} \supp(\varphi^{J_{i}^{\varphi}(r)})$ for $i = 1, 2, \ldots, p$.
    Cells are merged to form a parent-level cell by computing and factorizing the interactions between a cell and its corresponding proxy boundary, and the interactions between a cell and its adjacent cells
  }
  \label{fig:proxy}
\end{figure}
\begin{figure}[tb]
  \centering
  \includegraphics[width=0.30\linewidth]{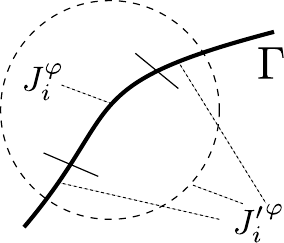}
  \caption{
      Correspondence between index set ${J_i^{\prime}}^{\varphi}$ and boundaries.
      The circle represented by the dashed line indicates a proxy boundary.
      The index ${J_i^{\prime}}^{\varphi}$ corresponds to the union of indices on the proxy boundary enclosing the cell associated with $J_i^{\varphi}$, and indices on the enclosed part of the adjacent cells of $J_i^{\varphi}$
  }
  \label{fig:proxy_j_prime}
\end{figure}
In this work, the subsequences $\tilde{J}_{i}^{\varphi}$ and $\tilde{J}_{i}^{\psi}$ and
coefficients $L_i$ and $R_i$ are computed by the proxy method \cite{MARTINSSON20051}
and interpolative factorization \cite{cheng2005compression}.
First, we prepare the notations related to the proxy method.
For details of the interpolative factorization, see Appendix \ref{ap:id}.
Figure \ref{fig:proxy} provides an overview of the proxy method.
The proxy method is the low-rank approximation method based on replacing far interactions with local interactions,
where far and adjacent interactions correspond to the off-diagonal blocks of \eqref{eq:block_linear}.
We can represent the part of $\Gamma$ with respect to $J_i^{\varphi}$ as $\bigcup_{r = 1}^{n} \supp(\varphi^{J_i^{\varphi} (r)})$
by the support of bases $\{ \varphi^r \}_{r = 1}^{n}$ on the boundary $\Gamma$.
We define the proxy boundary as a virtual boundary enclosing each $\bigcup_{r = 1}^{n} \supp(\varphi^{J_i^{\varphi} (r)})$,
which are the dashed lines in Figure \ref{fig:proxy}.
We can also approximate each proxy boundary with a polygon.
We define the piecewise linear basis on this polygon.
Let ${J_i^{\prime}}^{\varphi}$ be the sum of indices for the set of piecewise linear bases
on this approximated proxy boundary with respect to the $i$-th cell $J_i^{\varphi}$
and on the interior part of adjacent cells enclosed by the same proxy boundary as illustrated in Figure \ref{fig:proxy_j_prime}.
Similarly, the set of indices ${J_i^{\prime}}^{\psi}$ of the piecewise constant bases
on a proxy boundary with respect to $i$-th cell $J_i^{\psi}$
are defined.
These indices define the pair ${J_i^{\prime}} = ({J_i^{\prime}}^{\varphi}, {J_i^{\prime}}^{\psi})$.
For $i = 1, 2, \ldots, p$,
the interactions $M(J_i^{\prime}, J_i)$ and $M(J_i, J_i^{\prime})$ are also defined as same way for $M(\tilde{J}_i, \tilde{J}_j)$.

\begin{table}[h]
  \caption{Target matrix for the column-pivoted QR factorization to obtain each $R_{i}^{j}$}
  \label{tb:v}%
  \begin{tabular}{@{}ll@{}}
    \toprule
    Resulting matrix & Target for column-pivoted QR factorization \\
    \midrule
    $R_{i}^{1}$ & $\mqty(M_{\mathrm{lh}}^{1} ({J_i}^{\prime}, J_i) \\ M_{\mathrm{lh}}^{3} ({J_i}^{\prime}, J_i))$ \\
    $R_{i}^{2}$ & $\mqty(M_{\mathrm{rh}}^{1} ({J_i}^{\prime}, J_i) \\ M_{\mathrm{rh}}^{3} ({J_i}^{\prime}, J_i))$   \\
    $R_{i}^{3}$ & $\mqty(M_{\mathrm{lh}}^{2} ({J_i}^{\prime}, J_i) \\ M_{\mathrm{lh}}^{4} ({J_i}^{\prime}, J_i))$   \\
    $R_{i}^{4}$ & $\mqty(M_{\mathrm{rh}}^{2} ({J_i}^{\prime}, J_i) \\ M_{\mathrm{rh}}^{4} ({J_i}^{\prime}, J_i))$   \\
    \botrule
  \end{tabular}
\end{table}
\begin{table}[h]
  \caption{Target matrix for the column-pivoted QR factorization to obtain each $L_{i}^{j}$}
  \label{tb:u}%
  \begin{tabular}{@{}ll@{}}
    \toprule
    Resulting matrix & Target for column-pivoted QR factorization \\
    \midrule
    $L_{i}^{1}$    & $\mqty(M_{\mathrm{uh}}^{1} (J_i, {J_i}^{\prime}) & M_{\mathrm{uh}}^{2} (J_i, {J_i}^{\prime}))^{H}$ \\
    $L_{i}^{2}$    & $\mqty(M_{\mathrm{dh}}^{1} (J_i, {J_i}^{\prime}) & M_{\mathrm{dh}}^{2} (J_i, {J_i}^{\prime}))^{H}$   \\
    $L_{i}^{3}$    & $\mqty(M_{\mathrm{uh}}^{3} (J_i, {J_i}^{\prime}) & M_{\mathrm{uh}}^{4} (J_i, {J_i}^{\prime}))^{H}$   \\
    $L_{i}^{4}$    & $\mqty(M_{\mathrm{dh}}^{3} (J_i, {J_i}^{\prime}) & M_{\mathrm{dh}}^{4} (J_i, {J_i}^{\prime}))^{H}$   \\
    \botrule
  \end{tabular}
\end{table}
We now describe the calculation procedure for the left and right coefficients $L_i$ and $R_i$
and corresponding subsequences $\tilde{J}_i = (\tilde{J}_i^{\varphi}, \tilde{J}_i^\psi)$ for $i = 1, 2, \ldots, p$.
In \eqref{eq:a=usv},
$L_i$ and $R_i$ are written vaguely,
but to be precise, they are block diagonal matrices including four sub-matrices, such as \begin{equation}
  L_i =
  \mqty(
  L_{i}^{1} & & & \\
  & L_{i}^{2} & & \\
  & & L_{i}^{3} & \\
  & & & L_{i}^{4}
  ), \quad
  R_i =
  \mqty(
  R_{i}^{1} & & & \\
  & R_{i}^{2} & & \\
  & & R_{i}^{3} & \\
  & & & R_{i}^{4}
  ),
\end{equation}
where $L_{i}^{j} \in \mathbb{C}^{n \times k}$ and $R_{i}^{j} \in \mathbb{C}^{k \times n}$ for $j = 1, 2, 3, 4$.
In this work, we use column-pivoted QR factorization to compute the interpolative factorization.
  The interpolative factorization is a technique that decomposes a low-rank approximable matrix by using a selection of columns (or rows) to interpolate the entire matrix.
  In this study, it is used to obtain the skeleton matrix $S_{ij}$ $(i \neq j)$ as a submatrix of $A_{ij}$.
Tables \ref{tb:v} and \ref{tb:u} describe the target matrix for column-pivoted QR factorization to obtain each $R_{i}^{j}$ and $L_{i}^{j}$, respectively, for $j = 1, 2, 3, 4$.
In table \ref{tb:u}, $B^H$ stands for the complex conjugate and transpose of a matrix $B$.
The subsequences $\tilde{J}_{i}^{\varphi}$ and $\tilde{J}_{i}^{\psi}$ 
can be determined through the indices corresponding to the selected first $k$ pivots of the column-pivoted QR factorization,
because pivots are selected in descending order of the column norm of the target interaction matrix.
For example, in the PMCHWT formulation, $R_i^1$ and $R_i^2$ correspond to the same indices $\tilde{J}_{i}^{\psi}$.
Similarly, $R_i^3$ and $R_i^4$ correspond to the same indices $\tilde{J}_{i}^{\varphi}$ in the PMCHWT formulation.
In this work, we first compute $R_i^1$, $R_i^3$ and $\tilde{J}_{i}$ for the cell indices $J_i$,
then we compute $R_i^2$ and $R_i^4$ using the computed indices $\tilde{J}_{i}$.
More precisely, we use the permutation matrix $P$
that can be determined from the pivots corresponding to the computed indices $\tilde{J}_{i}$.
Instead of performing column-pivoted QR factorization on each target matrix $M$ to obtain $R_i^2$ and $R_i^4$,
we perform standard QR factorization on the permuted matrix $MP$,
then cut off the rank of the factorized matrix to $k$ in the same way as in the interpolative factorization.
A similar strategy is applied to compute $L_i^p$ for $p = 1, 2, 3, 4$.

The above method allows us to compute the shared left and right coefficients
needed for low-rank approximations of the off-diagonal parts of the linear equation \eqref{eq:block_linear}
in transmission problems for elastic waves.
The difference from the Neumann problem \cite{MATSUMOTO2025106148} is
that two types of bases $\{ \varphi^r \}$ and $\{ \psi^r \}$
are used for discretization, piecewise linear and piecewise constant,
and therefore the two index sets $\{ \tilde{J}_i^\varphi \}$ and $\{ \tilde{J}_i^\psi \}$
for low-rank approximations are required.

\subsection{Fast direct solver}
By using the low-rank approximation method described in the previous section,
we can apply almost the same fast direct solver formulation as that in \cite{MARTINSSON20051} and \cite{MATSUMOTO2025106148},
even when employing the versatile discretization based on mixed piecewise linear and constant bases.
The multi-level algorithm can be formulated by 
a straightforward extension of the rearrangement of matrices and vectors accompanying level transitions 
in the case of the two-block \cite{MATSUMOTO2025106148} to the four-block case.

The single-level algorithm for the proposed fast direct solver is denoted as follows.
The block linear system \eqref{eq:block_linear} becomes
\begin{equation}
  \mqty(
  A_{11} & L_{1} S_{12} R_{2} & \cdots & L_{1} S_{1p} R_{p} \\
  L_{2} S_{21} R_{1} & A_{22} & \cdots & L_{2} S_{2p} R_{p} \\
  \vdots & \vdots &\ddots &\vdots \\
  L_{p} S_{p1} R_{1} & L_{p} S_{p2} R_{2} & \cdots & A_{pp}
  )
  \mqty(
  x_1 \\
  x_2 \\
  \vdots \\
  x_p
  )
  =
  \mqty(
  f_1 \\
  f_2 \\
  \vdots \\
  f_p
  ),
  \label{eq:low-rank}
\end{equation}
by the low-rank approximations of off-diagonal blocks with the shared $L_i$ and $R_i$.
In this work, the shared $L_i$ and $R_i$ are computed by the proxy method described in the previous section,
but we can use other low-rank approximation techniques to compress the system to a smaller system
if it is a weak admissibility type approximation.
Note that the size of the coefficient matrix for system \eqref{eq:low-rank} is $4np \times 4np$,
where $n = |J_i^{\varphi}| = |J_i^{\psi}|$ and $p$ is the number of cells in the leaf level; therefore, $np = N$.
By multiplying $\diag(R_{1}, R_{2}, \ldots, R_{p})$
and $\diag(A_{11}^{-1}, A_{22}^{-1}, \ldots, A_{pp}^{-1})$ from the left, we have
\begin{multline}
  \mqty(
  R_{1}x_1 \\
  R_{2}x_2 \\
  \vdots \\
  R_{p}x_p
  )
  +
  \mqty(
  R_{1}A_{11}^{-1}L_{1} & & & \\
   & R_{2}A_{22}^{-1}L_{2} & & \\
   & & \ddots & \\
   & & & R_{p}A_{pp}^{-1}L_{p}
  )
  \mqty(
  0 & S_{12} R_{2} & \cdots & S_{1p} R_{p} \\
  S_{21} R_{1} & 0 & \cdots & S_{2p} R_{p} \\
  \vdots & \vdots &\ddots &\vdots \\
  S_{p1} R_{1} & S_{p2} R_{2} & \cdots & 0
  )
  \mqty(
  x_1 \\
  x_2 \\
  \vdots \\
  x_p
  )
  \\
  =
  \mqty(
  R_{1}A_{11}^{-1}f_1 \\
  R_{2}A_{22}^{-1}f_2 \\
  \vdots \\
  R_{p}A_{pp}^{-1}f_p
  ).
\end{multline}
Let $y_i := R_i x_i$ and $\tilde{A}_{i} := \qty(R_i A_{ii}^{-1} L_i)^{-1}$.
By multiplying $\diag(\tilde{A}_{1}, \tilde{A}_{2}, \ldots, \tilde{A}_{p})$ from the left,
we obtain compressed block linear equations expressed as
\begin{equation}
  \mqty(
  \tilde{A}_{1} & S_{12} & \cdots & S_{1p} \\
  S_{21} & \tilde{A}_{2} & \cdots & S_{2p} \\
  \vdots & \vdots &\ddots &\vdots \\
  S_{p1} & S_{p2} & \cdots & \tilde{A}_{p}
  )
  \mqty(
  y_1 \\
  y_2 \\
  \vdots \\
  y_p
  )
  =
  \mqty(
  \tilde{A}_{1}R_{1}A_{11}^{-1}f_1 \\
  \tilde{A}_{2}R_{2}A_{22}^{-1}f_2 \\
  \vdots \\
  \tilde{A}_{p}R_{p}A_{pp}^{-1}f_p
  ),
  \label{eq:compressed}
\end{equation}
where the vectors $\{y_i \}_{i = 1}^{p}$ are dealt with as new unknowns.
Because the size of the coefficient matrix for this system is $4kp \times 4kp$ ($k \ll n$),
this system can be solved faster than the original one \eqref{eq:block_linear} by a linear algebra solver package such as ZGESV of Lapack.

After solving the above equation for $\{y_i \}_{i = 1}^{p}$,
we obtain the original solution vectors $\{x_i \}_{i = 1}^{p}$ with the following calculations.
Multiplying $\diag(A_{11}^{-1} L_1, A_{22}^{-1} L_i, \ldots, A_{pp}^{-1} L_p)$ from the left of \eqref{eq:compressed} yields
\begin{equation}
  \mqty(
  A_{11}^{-1} L_1 &  &  & \\
   & A_{22}^{-1} L_2 &  & \\
  &  &\ddots & \\
   &  & & A_{pp}^{-1} L_p
  )
  \mqty(
  \tilde{A}_{1} & S_{12} & \cdots & S_{1p} \\
  S_{21} & \tilde{A}_{2} & \cdots & S_{2p} \\
  \vdots & \vdots &\ddots &\vdots \\
  S_{p1} & S_{p2} & \cdots & \tilde{A}_{p}
  )
  \mqty(
  y_1 \\
  y_2 \\
  \vdots \\
  y_p
  )
  =
  \mqty(
  A_{11}^{-1} L_1 \tilde{A}_{1}R_{1}A_{11}^{-1}f_1 \\
  A_{22}^{-1} L_2 \tilde{A}_{2}R_{2}A_{22}^{-1}f_2 \\
  \vdots \\
  A_{pp}^{-1} L_p \tilde{A}_{p}R_{p}A_{pp}^{-1}f_p
  ).
  \label{eq:eliminate_target}
\end{equation}
We want to eliminate the terms with respect to $S_{ij}$ from this equation
by finding another expression that includes $\{ x_i \}_{i = 1}^{p}$ for
\begin{equation}
  \mqty(
  A_{11}^{-1} L_1 &  &  & \\
   & A_{22}^{-1} L_2 &  & \\
  &  &\ddots & \\
   &  & & A_{pp}^{-1} L_p
  )
  \mqty(
   & S_{12} & \cdots & S_{1p} \\
  S_{21} &  & \cdots & S_{2p} \\
  \vdots & \vdots &\ddots &\vdots \\
  S_{p1} & S_{p2} & \cdots & 
  )
  \mqty(
  y_1 \\
  y_2 \\
  \vdots \\
  y_p
  ).
\end{equation}
We therefore multiply both sides of \eqref{eq:low-rank}
by $\diag(A_{11}^{-1}, A_{22}^{-1}, \ldots, A_{pp}^{-1})$ from the left to obtain the relation
\begin{equation}
  \mqty(
  A_{11}^{-1} L_1 &  &  & \\
   & A_{22}^{-1} L_2 &  & \\
  &  &\ddots & \\
   &  & & A_{pp}^{-1} L_p
  )
  \mqty(
   & S_{12} & \cdots & S_{1p} \\
  S_{21} &  & \cdots & S_{2p} \\
  \vdots & \vdots &\ddots &\vdots \\
  S_{p1} & S_{p2} & \cdots & 
  )
    \mqty(
  y_1 \\
  y_2 \\
  \vdots \\
  y_p
  )
  =
  -
  \mqty(
  x_1 \\
  x_2 \\
  \vdots \\
  x_p
  )
  +
  \mqty(
  A_{11}^{-1} f_1 \\
  A_{22}^{-1} f_2 \\
  \vdots \\
  A_{pp}^{-1} f_p
  ).
\end{equation}
Substituting this into \eqref{eq:eliminate_target}, we have
\begin{equation}
  \mqty(
  A_{11}^{-1} L_1 \tilde{A}_{1} y_1 \\
  A_{22}^{-1} L_2 \tilde{A}_{2} y_2 \\
  \vdots \\
  A_{pp}^{-1} L_p \tilde{A}_{p} y_p
  )
  -
  \mqty(
  x_1 \\
  x_2 \\
  \vdots \\
  x_p
  )
  +
  \mqty(
  A_{11}^{-1} f_1 \\
  A_{22}^{-1} f_2 \\
  \vdots \\
  A_{pp}^{-1} f_p
  )
  =
  \mqty(
  A_{11}^{-1} L_1 \tilde{A}_{1}R_{1}A_{11}^{-1}f_1 \\
  A_{22}^{-1} L_2 \tilde{A}_{2}R_{2}A_{22}^{-1}f_2 \\
  \vdots \\
  A_{pp}^{-1} L_p \tilde{A}_{p}R_{p}A_{pp}^{-1}f_p
  ).
\end{equation}
By simplifying this, we obtain the conversion relation from
$\{ y_i \}_{i = 1}^{p}$ to $\{ x_i \}_{i = 1}^{p}$ as
\begin{equation}
  \mqty(
  x_1 \\
  x_2 \\
  \vdots \\
  x_p
  )
  =
  \mqty(
  A_{11}^{-1}f_1 \\
  A_{22}^{-1}f_2 \\
  \vdots \\
  A_{pp}^{-1}f_p
  )
  -
  \mqty(
  A_{11}^{-1}L_{1} \tilde{A}_{1}R_{1}A_{11}^{-1}f_1 \\
  A_{22}^{-1}L_{2}\tilde{A}_{2}R_{2}A_{22}^{-1}f_2 \\
  \vdots \\
  A_{pp}^{-1}L_{p}\tilde{A}_{p}R_{p}A_{pp}^{-1}f_p
  )
  +
  \mqty(
  A_{11}^{-1}L_{1} \tilde{A}_{1} y_1 \\
  A_{22}^{-1}L_{2} \tilde{A}_{2} y_2 \\
  \vdots \\
  A_{pp}^{-1}L_{p} \tilde{A}_{p} y_p
  ).
\end{equation}

\begin{remark}
When $N$ is the size of the system, 
the above single-level fast direct solver has $O( (4kp)^3) = O(N^3)$ complexity because the number of cells is $p = O(N)$.
Here, we assume the number of skeletons $k = O(1)$ and $k \ll n$ by a sufficient depth of tree partitioning.
If this assumption is accepted, we can formulate the multi-level fast direct solver that can compress \eqref{eq:block_linear} to a system of size $O(1)$
because the merged $S_{ij}$ in \eqref{eq:compressed} can be low-rank approximated recursively.
This results in an $O(N)$ fast direct solver, which is used in the following numerical examples.
For details of the multi-level algorithm, see \cite{MATSUMOTO2025106148}.
\end{remark}

\section{Numerical examples} \label{sec:Numerical}
The examples in this section demonstrate the performance of the proposed fast direct solver based on the proxy method in Section \ref{sec:proposed}
for the transmission problems for elastic waves.
This fast direct solver has versatility, that is, the same program codes can be used for inclusions that have
a smooth boundary or a boundary with corner points.

\begin{figure}[tb]
  \centering
  \includegraphics[width=0.4\linewidth]{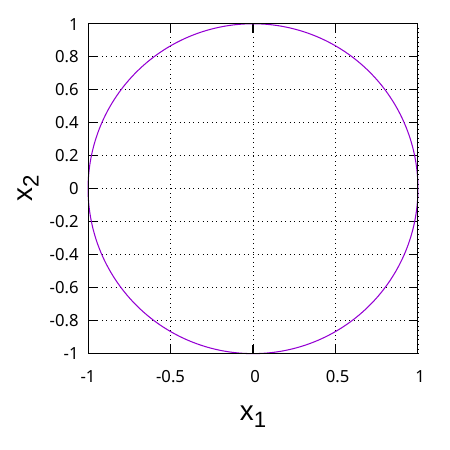}
  \caption{Geometry for the unit circle}
  \label{fig:mesh_circle}
\end{figure}

\begin{figure}[tb]
  \centering
  \includegraphics[width=0.4\linewidth]{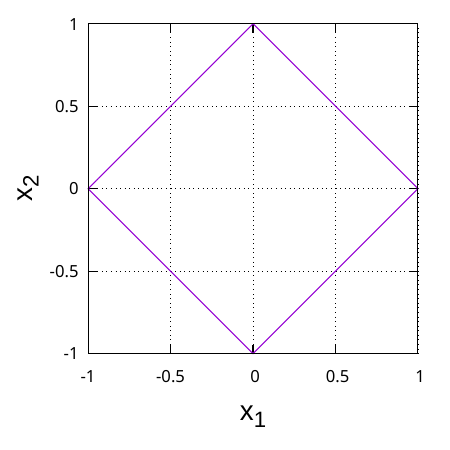}
  \caption{Geometry for the square}
  \label{fig:mesh_square}
\end{figure}

In what follows, a longitudinal plane wave is used as an incident wave for all numerical examples.
The performance of the fast direct solver for the formulations of the boundary element method of the PMCHWT formulation \eqref{eq:gal_pm1} and \eqref{eq:gal_pm2}
and the Burton--Miller formulation \eqref{eq:gal_bm1} and \eqref{eq:gal_bm2} is demonstrated.
The DOF is equal to $4 N$,
where $N = N_n = N_e$ and $N_n$ and $N_e$ are the number of nodes and elements on the approximated $\Gamma$, respectively.
We call the corresponding conventional boundary element method ``Conv'' in this paper.
The choice of basis and test functions used for Conv is identical to that used for the fast direct solver.
The system of linear equations derived from Conv is solved using a standard LU factorization with partial pivoting.
In this numerical example, a unit circle and a square are used for an inclusion, as indicated in Figures \ref{fig:mesh_circle} and \ref{fig:mesh_square}.
The multi-level algorithm for the proposed fast direct solver is used.
A binary tree is used for a tree partitioning of the indices of nodes and elements.
The bottom level of this binary tree is called the leaf level, which has $p$ cells for both the node tree and element tree.
The rank after low-rank approximation is expressed as the number of skeletons.
This number of skeletons increases with each level in the binary tree
by $1.15$ times that for the previous level in all numerical examples.
For example, if the number of skeletons at the leaf level is 30,
it increases to 34 at the level above the leaf level,
and then to 39 at the level above that.
We refer to the number of skeletons in the leaf level as the ``initial rank.'' This initial rank is used for the low-rank approximations when computing $L_i^{1}, L_i^{2}, L_i^{3}, L_i^{4}, R_i^{1}, R_i^{2}, R_i^{3}$, and $R_i^{4}$ $(i = 1, 2, \ldots, p)$ at the leaf level.
For computation of the numerical examples,
a single compute node of the supercomputer TSUBAME4.0 at the Institute of Science Tokyo,
which has two 96-core AMD EPYC 9654 CPUs (total number of cores: 192;
maximum frequency: 3708 MHz) and 768 GB of 0.92-TB/s DDR5 memory,
was used.
The implementation was parallelized by OpenMP.
We utilized all 192 CPU cores for all numerical examples.
The fast direct solver part was implemented in C++,
and the calculation of the discretized layer potential was implemented in Fortran.
The compilers used were g++ and gfortran, respectively.

\subsection{Verification of implementation of the fast direct solver} \label{sec:veri_fds}
\begin{figure}[tb]
  \centering
  \includegraphics[width=0.8\linewidth]{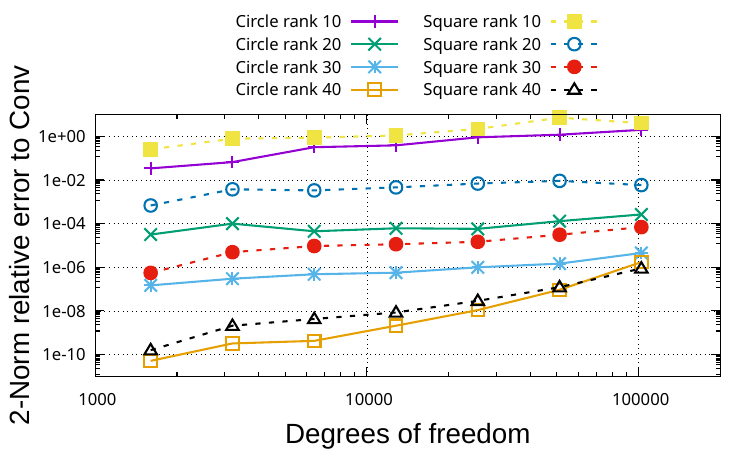}
  \caption{Two-norm relative error defined by \eqref{eq:relative_norm} for the fast direct solver based on the PMCHWT formulation.
    In this figure, ``Circle'' and ``Square'' correspond to the geometry.
    Furthermore, the rank numbers (e.g., ``rank 10'') stand for the number of skeletons for $L_i^{1}, L_i^{2}, L_i^{3}, L_i^{4}, R_i^{1}, R_i^{2}, R_i^{3}$, and $R_i^{4}$ $(i = 1, 2, \ldots, p)$ on the leaf level
  }
  \label{fig:error_pmchwt_fds}
\end{figure}
\begin{figure}[tb]
  \centering
  \includegraphics[width=0.8\linewidth]{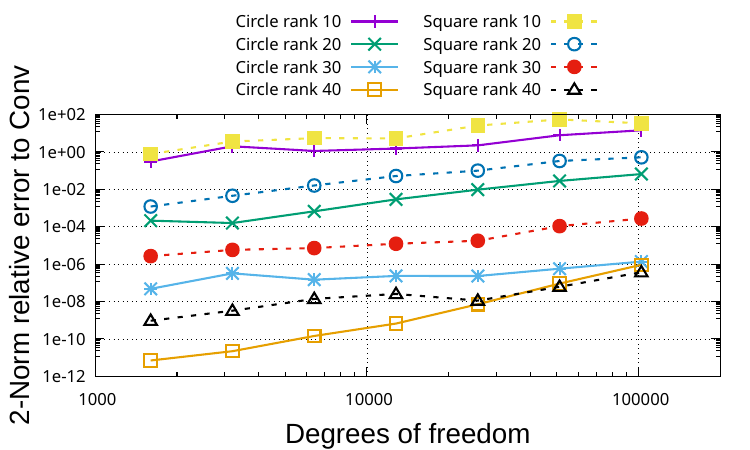}
  \caption{Two-norm relative error defined by \eqref{eq:relative_norm} for the fast direct solver based on the Burton--Miller formulation. The legend text has the same meaning as described in the Figure \ref{fig:error_pmchwt_fds} caption}
  \label{fig:error_bm_fds}
\end{figure}
At first, implementation of the fast direct solver is verified by
comparing the results to the numerical solutions of the conventional boundary element method (Conv).
The implementation of Conv is verified in Appendix \ref{ap:analytical}.

The parameters used in this verification are as follows.
The frequency of the incident wave is $\omega = 4$.
The wave speeds in $\Omega_{0}$ are $c_{L}^{(0)} = \sqrt{3}$, $c_{T}^{(0)} = 1$.
The wave speeds in $\Omega_{1}$ are $c_{L}^{(1)} = 3$, $c_{T}^{(1)} = 1.5$.
The densities of $\Omega_{0}$ and $\Omega_{1}$ are $\rho^{(0)} = 1$ and $\rho^{(1)} = 2$, respectively.

The numerical solution vector $x$ is expressed as
\[
x := ( u_1^1, u_1^2, \ldots, u_1^{N_n}, u_2^1, u_2^2, \ldots, u_2^{N_n}, t_1^1, t_1^2, \ldots, t_1^{N_e}, t_2^1, t_2^2, \ldots, t_2^{N_e} ) \in \mathbb{C}^{4N},
\]
on $\Gamma$.
The notations $x^{\mathrm{Conv}}$ and $x^{\mathrm{FDS}}$ indicate the numerical solutions computed by the conventional boundary element method and by the fast direct solver, respectively.
Figures \ref{fig:error_pmchwt_fds} and \ref{fig:error_bm_fds} show the accuracy of the fast direct solver compared to that of Conv, as calculated by
\begin{equation}
  \frac{\norm{x^{\mathrm{FDS}} - x^{\mathrm{Conv}}}_2}{\norm{x^{\mathrm{Conv}}}_2}
  \label{eq:relative_norm}
\end{equation}
for the PMCHWT and Burton--Miller formulations, respectively.
Here, $\norm{\cdot}_2$ means the $2$-norm of a vector.

Figures \ref{fig:error_pmchwt_fds} and \ref{fig:error_bm_fds} demonstrate the above relative $2$-norm
up to $102,400$ DOFs.
Here, this limitation on the DOF arises from the memory capacity used in Conv.
Figures \ref{fig:error_pmchwt_fds} and \ref{fig:error_bm_fds} show that
the numerical solutions for the fast direct solvers based on the PMCHWT and Burton--Miller formulations well match
that of Conv when the initial rank is $30$ or $40$.
In the case of a initial rank of 30 or 40, almost all relative $2$-norm errors are smaller than $1.0 \times 10^{-4}$.
However, there is a deterioration in accuracy when the initial rank is $10$ or $20$.
This indicates that the accuracy of the fast direct solver
can be controlled by the number of skeletons in low-rank approximations.

\begin{remark}
A decrease in accuracy of the fast direct solver against Conv as the DOF
increases has also been observed in cavity scattering \cite{MATSUMOTO2025106148},
and therefore is not unique to transmission problems.
A similar deterioration of accuracy of the fast direct solvers based on the proxy method is also observed in the Helmholtz transmission problems discretized by the collocation method with piecewise constant bases \cite{YasuhiroMATSUMOTO202308-231124}.
In the result of the Helmholtz transmission problems, the solution of the fast direct solver is compared to the analytical solution on a unit circle.
However, in our experience, when the Helmholtz transmission problem is discretized using
a low-order Nystr\"om method, the accuracy degradation issue due to the increase in the DOF does not occur.
The cause might be related to the discretization method, but presently it remains unidentified.
It should also be noted that elastic and acoustic scattering are different problems.
\end{remark}

\subsection{Time and memory complexities of the proposed solver} \label{sec:complexity}
This section describes the performance of the proposed fast direct solver, which has $O(N)$ complexity at fixed (low) frequencies.
The frequency of the incident wave $\omega$, 
the wave speeds $c_{L}^{(m)}$ and $c_{T}^{(m)}$,
and the densities $\rho^{(m)}$ for $m = 0, 1$ are the same as those given in the previous section.
Considering the results for the accuracy for the fast direct solvers in the previous section,
the performance of the proposed fast direct solver is measured for the cases of initial ranks 30 and 40.

\begin{figure}[tb]
  \centering
  \includegraphics[width=0.8\linewidth]{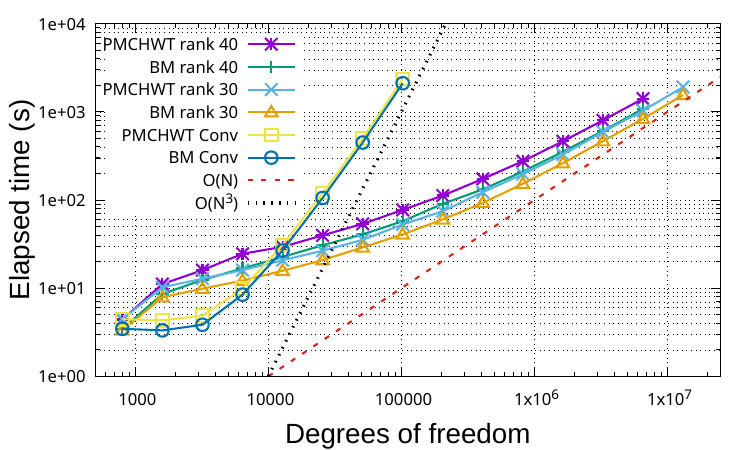}
  \caption{Computational time for proposed fast direct solvers for a circular inclusion.
    In the legend, ``PMCHWT'' and ``BM'' mean the PMCHWT and Burton--Miller formulations, respectively, and
    ``rank 30'' and ``rank 40'' correspond to initial ranks 30 and 40, respectively.
    Furthermore, ``Conv'' corresponds to the conventional boundary element method solved by a standard LU factorization
  }
  \label{fig:time_nlev_circle}
\end{figure}
\begin{figure}[tb]
  \centering
  \includegraphics[width=0.8\linewidth]{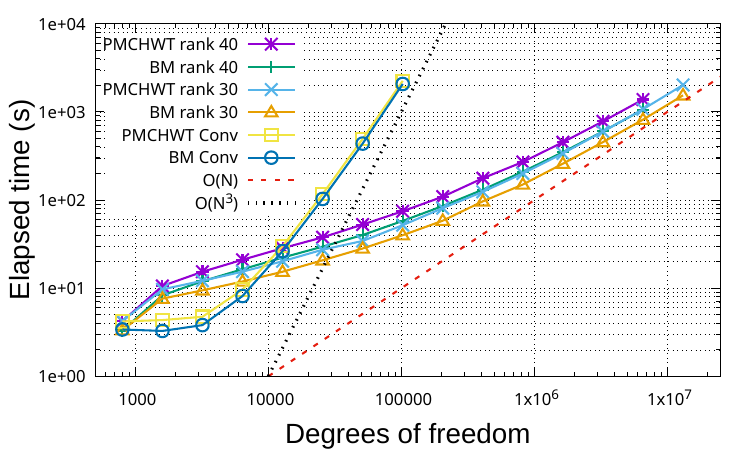}
  \caption{Computational time for proposed fast direct solvers for a square inclusion.
    The legend entries have the same meanings as those given in the caption of Figure \ref{fig:time_nlev_circle}
  }
  \label{fig:time_nlev_square}
\end{figure}
\begin{figure}[tb]
  \centering
  \includegraphics[width=0.9\linewidth]{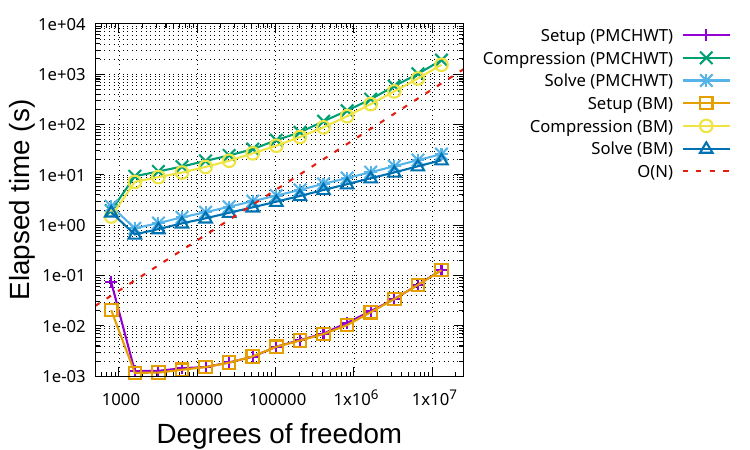}
  \caption{Computational time of each step for proposed fast direct solvers for a circular inclusion.
    This figure shows the case of ``rank 30'' of both PMCHWT and Burton--Miller (BM) formulations of Figure \ref{fig:time_nlev_circle}.
    In the legend, the step of ``Setup'' corresponds to the preprocessing stage, including mesh generation and indexing.
    The ``Compress'' step corresponds to the recursive compression of the system linear equations using the proxy surface method. The right-hand side at the leaf level is also computed in the ``Compress'' step.
    The ``Solve'' step involves solving the compressed system, which includes computing the remaining off-diagonal blocks and reconstructing the original solution from the compressed one
  }
  \label{fig:time_step_nlev_circle}
\end{figure}
\begin{figure}[tb]
  \centering
  \includegraphics[width=0.9\linewidth]{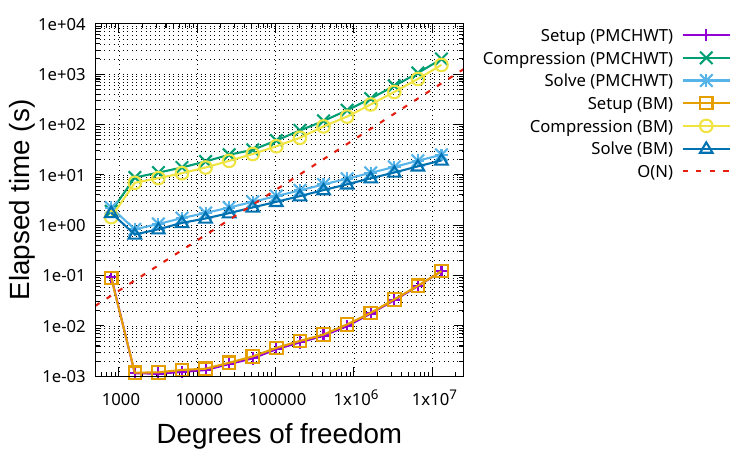}
  \caption{Computational time of each step for proposed fast direct solvers for a square inclusion.
    This figure shows the case of ``rank 30'' of both PMCHWT and Burton--Miller (BM) formulations of Figure \ref{fig:time_nlev_square}.
    The legend entries have the same meanings as those given in the caption of Figure \ref{fig:time_step_nlev_circle}
  }
  \label{fig:time_step_nlev_square}
\end{figure}
Figures \ref{fig:time_nlev_circle} and \ref{fig:time_nlev_square}
show the computational time for the circle and square geometries for the proposed fast direct solvers.
In these figures, the elapsed time includes all computations from mesh generation to obtaining the numerical solution.
Figure \ref{fig:time_nlev_circle} shows that the fast direct solvers based on both formulations
have a computational complexity of $O(N)$.
Up to about 10,000 DOFs, Conv is faster than the fast direct solver in both the PMCHWT and Burton--Miller formulations
because the speedup from parallelization is greater in Conv, which has a simpler algorithm.
However, because Conv requires $O(N^2)$ memory and $O(N^3)$ floating point operations,
Conv cannot solve large-scale problems with more than
102,400 DOFs due to insufficient memory capacity.
Also, the computational time needed to solve the 102,400
DOF problem with Conv is longer than
that to solve the 13,107,200 DOF problem
with the fast direct solver in the case of initial rank 30.
The results in Figure \ref{fig:time_nlev_square} are almost the same as those in Figure \ref{fig:time_nlev_circle}.
This implies that the inclusion shape does not affect the computational time in the proposed fast direct solver.

  Figures \ref{fig:time_step_nlev_circle} and \ref{fig:time_step_nlev_square}
  show the computational time of each step for the circle and square geometries for the proposed fast direct solvers for both PMCHWT and Burton--Miller formulations in the case of ``rank 30'', respectively.
  These figures indicate that the dominant part of the computational time is the recursive compression of the system of linear equations based on the proxy method.
  The ``Solve'' step exhibits a relatively long computational time compared to ``$\mathrm{Time}^{\text{2nd--10th}}$'' in Table \ref{tb:second_and_after} (presented later).
  This is because the initial execution of this step involves the computation of interactions corresponding to the remaining off-diagonal blocks of the compressed system.
  Regarding the time complexity, as noted at the beginning of Section \ref{sec:Numerical}, the size of the compressed system in the proposed solver increases with the tree depth.
  This could potentially violate the $O(N)$ time complexity as the mesh is refined. However, as demonstrated in these figures, the ``Solve'' step is consistently faster than the ``Compress'' step. At least within the scope of our numerical experiments, we did not observe clear evidence that this effect leads to a breakdown of the $O(N)$ scaling upon mesh refinement.

\begin{figure}[tb]
  \centering
  \includegraphics[width=0.8\linewidth]{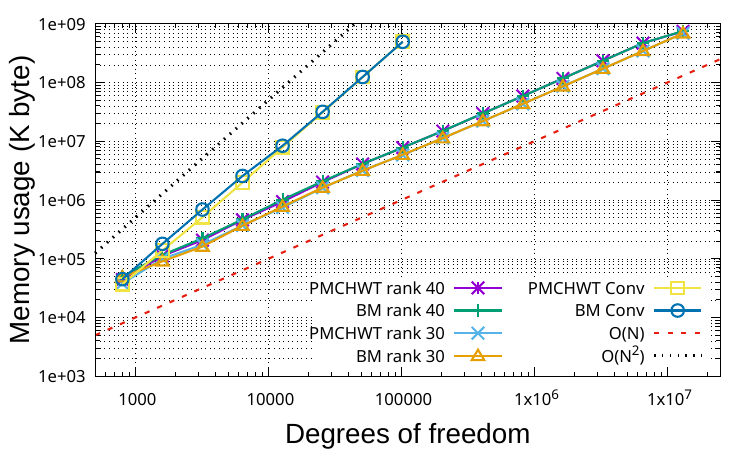}
  \caption{Memory usage for proposed fast direct solvers for a circular inclusion.
      The legend entries have the same meanings as those given in the caption of Figure \ref{fig:time_nlev_circle}
  }
  \label{fig:memory_nlev_circle}
\end{figure}
\begin{figure}[tb]
  \centering
  \includegraphics[width=0.8\linewidth]{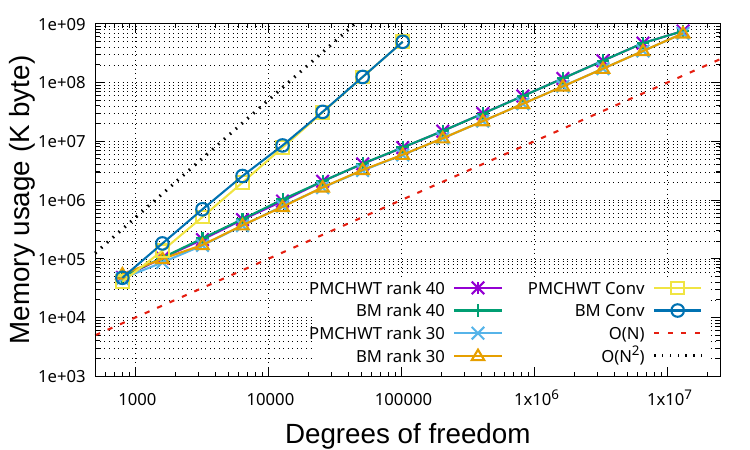}
  \caption{Memory usage for proposed fast direct solvers for a square inclusion.
      The legend entries have the same meanings as those given in the caption of Figure \ref{fig:time_nlev_square}
  }
  \label{fig:memory_nlev_square}
\end{figure}
Figures \ref{fig:memory_nlev_circle} and \ref{fig:memory_nlev_square} show the memory usage of the proposed fast direct solvers for the circle and square geometries, respectively. Memory usage was investigated using the {\it maximum resident set size} obtained from the GNU time command.
From these figures, we see that the proposed fast direct solvers for both the PMCHWT and Burton-Miller formulations demonstrate $O(N)$ memory usage as the mesh is refined, whereas the memory usage of Conv for both formulations is $O(N^2)$. Furthermore, similar to the computational time, changes in the shape of the scatterer have almost no effect on the memory usage.

\begin{table}[h]
  \caption{Average ratios of elapsed time obtained by the Burton--Miller (BM) formulation to that by the PMCHWT formulation}
  \label{tb:ratio}
  \begin{tabular}{@{}lll@{}}
    \toprule
    Inclusion shape & Cases & Average ratio ($\mathrm{Time}^{\mathrm{BM}}$/$\mathrm{Time}^{\mathrm{PMCHWT}}$) \\
    \midrule
    Circle & Initial rank 30 & 0.782 \\
    & Initial rank 40 & 0.762 \\
    & (Conv) & (0.842) \\
    Square & Initial rank 30 & 0.769 \\
    & Initial rank 40 & 0.775 \\
    & (Conv) & (0.847) \\
    \botrule
  \end{tabular}
\end{table}
Furthermore, the fast direct solver based on the Burton--Miller formulation outperforms that based on the PMCHWT formulation.
More precisely, the computational time for the former
is approximately $20\%$ shorter than that for the latter,
as summarized in Table \ref{tb:ratio}.

The fast direct solver is characterized by
its efficiency in handling multiple right-hand sides of the equations after the initial solution.
To evaluate this performance, we measured the total computational time
for 10 right-hand sides under the same conditions.
This scenario included an initial incident wave followed by nine additional waves
from angles different than the initial one.
Table \ref{tb:second_and_after} shows the measurement results when the shape of the inclusion is square.
From Table \ref{tb:second_and_after}, we see that the proposed fast direct solvers can efficiently solve the multiple right-hand side problems.
\begin{table}[h]
  \caption{
    Computational times for fast direct solver based on the Burton--Miller (BM) and PMCHWT formulations for 10 right-hand sides against to the DOFs.
    Conv stands for the computational times for the conventional boundary element method based on the BM and PMCHWT formulations, which are the references.
    The proposed fast direct solvers are the cases of initial ranks 30 and 40.
    We define $\mathrm{Time^{1st}}$ as the elapsed time to solve with respect to the first right-hand side.
    $\mathrm{Time}^{\text{2nd--10th}}$ represents the total time to solve with respect to the 2nd through 10th right-hand sides.
    The time units are seconds in this table.
    In our parallelization with multiple right-hand side problems,
    the second right-hand side onward was implemented by allocating a CPU core per equation right-hand side,
    which results in idle CPU cores.
  }
  \label{tb:second_and_after}
  \begin{tabular}{@{}llllllll@{}}
    \toprule
    & & \multicolumn{2}{l}{DOF $=12,800$} & \multicolumn{2}{l}{DOF $=102,400$} & \multicolumn{2}{l}{DOF $=1,638,400$} \\
    Formulation & Cases & $\mathrm{Time^{1st}}$ & $\mathrm{Time}^{\text{2nd--10th}}$ &  $\mathrm{Time^{1st}}$ & $\mathrm{Time}^{\text{2nd--10th}}$ & $\mathrm{Time^{1st}}$ & $\mathrm{Time}^{\text{2nd--10th}}$  \\
    \midrule
    BM & Initial rank 30 & 15.4 & 0.0250 & 39.9 & 0.217 & 259 & 3.39 \\
    & Initial rank 40 & 22.2 & 0.0370 & 58.0 & 0.295 & 349 & 4.77 \\
    & (Conv) & (26.3) & (0.151) & (2070) & (7.46) & - & - \\
    PMCHWT & Initial rank 30 & 20.2 & 0.0197 & 51.9 & 0.219 & 338 & 3.49 \\
    & Initial rank 40 & 28.4 & 0.0350 & 74.7 & 0.308 & 453 & 4.87 \\
    & (Conv) & (30.0) & (0.157) & (2250) & (6.60) & - & - \\
    \botrule
  \end{tabular}
\end{table}

\subsection{Performance for various calculation conditions} \label{sec:parameter}
In this section, we show that the computational time for the proposed fast direct solver is relatively robust against changes in calculation conditions, such as the density of the inclusion.
This is also a feature of fast direct solvers,
while the computational time for iterative solvers is generally highly dependent on the calculation conditions without preconditioning.
Let the frequency be $\omega = 4, 8$,
and let the density of the inclusion be $\rho^{(1)} = 0.1, 0.2, 0.3, \ldots, 9.9, 10.0$ for each frequency.
The density of the background medium is set to $\rho^{(0)} = 1$.
  The DOF is set to $12,800$ ($N = N_n = N_e = 3200$) for $\omega = 4$ and $25,600$ ($N = N_n = N_e = 6400$) for $\omega = 8$, respectively.
The wave speeds $c_{L}^{(j)}$ and $c_{T}^{(j)}$ in $\Omega_{j}$ $(j = 1, 2)$ remain the same as those used in Section \ref{sec:veri_fds}.
Considering the accuracy results for the fast direct solvers discussed in the previous section,
the initial rank of the fast direct solver should be set to 30 to
balance the speed and accuracy.

\begin{figure}[tb]
  \centering
  \includegraphics[width=1.0\linewidth]{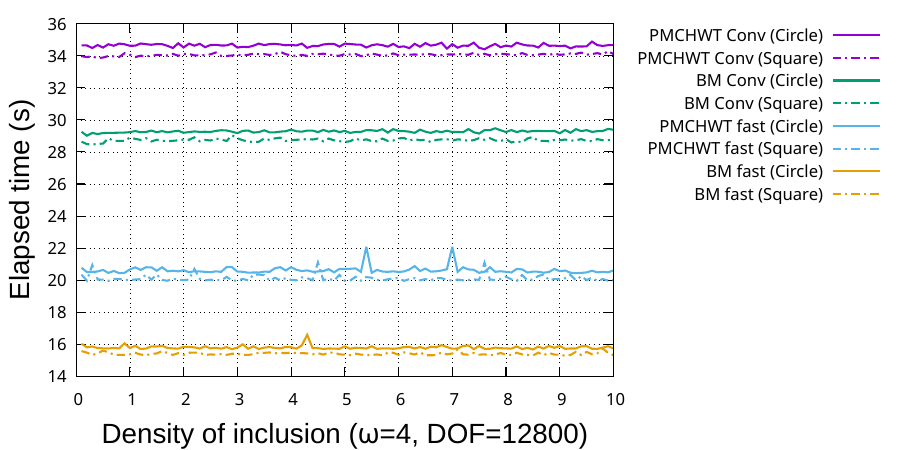}
  \caption{Computational time for DOF=12,800 with frequency $\omega = 4$ for various densities}
  \label{fig:time_rho_o4}
\end{figure}

\begin{figure}[tb]
  \centering
  \includegraphics[width=1.0\linewidth]{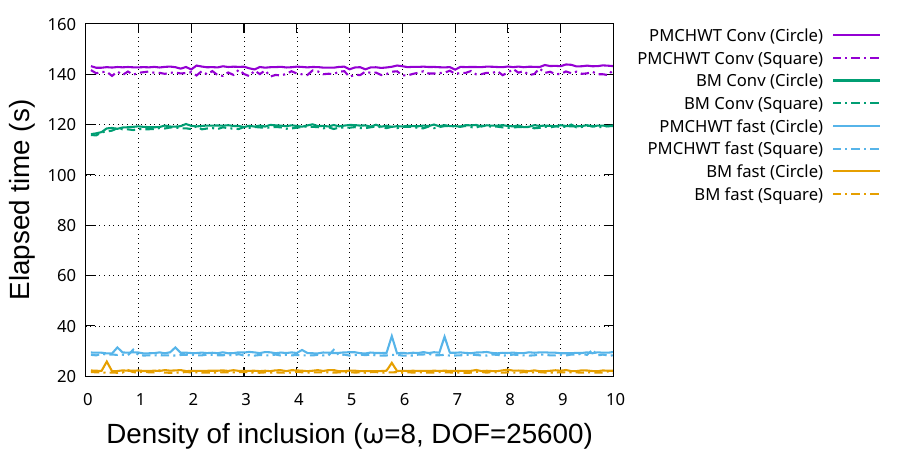}
  \caption{Computational time for DOF=25,600 with frequency $\omega = 8$ for various densities}
  \label{fig:time_rho_o8}
\end{figure}

Figures \ref{fig:time_rho_o4} and \ref{fig:time_rho_o8} show
the computational time for the proposed fast direct solvers
for $\omega =4$ and $\omega=8$, respectively.
The results for both circular and square inclusions are shown in both figures.
These figures show that the computational time for the proposed fast direct solvers is relatively stable against changing parameters.
Furthermore, These figures also show that, similar to the results given in the previous section, whether the inclusion shape is circular or square has little effect on the computational time.

\begin{figure}[tb]
  \centering
  \includegraphics[width=1.0\linewidth]{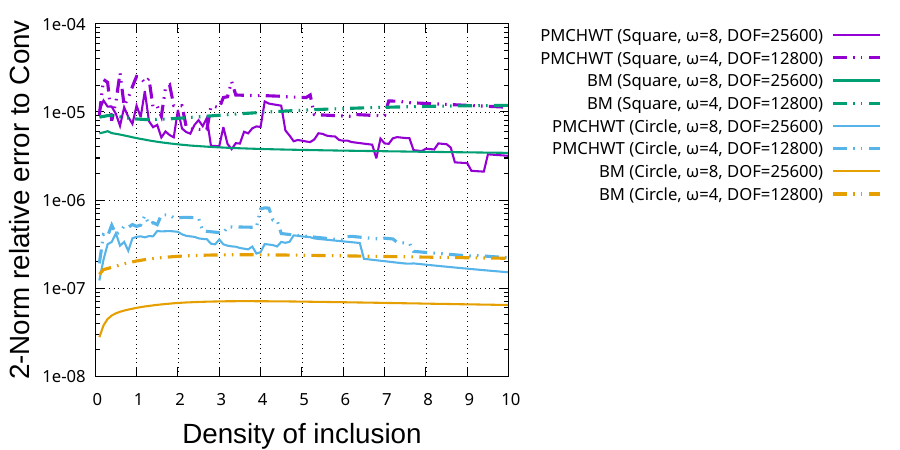}
  \caption{Relative 2-norm error \eqref{eq:relative_norm} for various density of inclusion. This figure corresponds to Figures \ref{fig:time_rho_o4} and \ref{fig:time_rho_o8}}
  \label{fig:error_rho}
\end{figure}

Figure \ref{fig:error_rho} shows
the relative 2-norm error \eqref{eq:relative_norm} against Conv,
which corresponds to the cases of 
Figures \ref{fig:time_rho_o4} and \ref{fig:time_rho_o8}.
The accuracy for the circular inclusion is higher than that for the square inclusion.
However, even for the square inclusion, an accuracy of $1.0 \times 10^{-4}$ is maintained.
The figures show that the accuracy is almost the same for the Burton--Miller and PMCHWT formulations,
suggesting that the low-rank approximation works well in both formulations.

\section{Conclusion and future work} \label{sec:Conclusion}
This work developed a fast direct solver for transmission problems for elastic waves in two dimensions.
The proposed fast direct solver is based on low-rank approximation of the off-diagonal blocks of the coefficient matrix for the discretized system via the proxy method.
The proposed fast direct solvers were formulated based on Burton--Miller type boundary integral equations and PMCHWT type boundary integral equations.
In these formulations, the displacement is approximated by piecewise linear bases
and the traction is approximated by piecewise constant bases.
This is a natural choice considering the continuity of displacement and discontinuity of stress on the boundary.
However, it prevents the use of the Calder\'on preconditioner, which is typical analytical preconditioning.
This therefore makes the application of the fast direct solvers appealing.
The proposed fast direct solver employs two independent index sets for the piecewise linear and constant basis to deal with this discretization strategy in the proxy method.
Numerical examples showed the correctness and performance of the proposed fast direct solvers based on the PMCHWT and Burton--Miller formulations.
For example, the solvers have $O(N)$ complexity at fixed low frequencies and are relatively robust against changing the density of an inclusion.
The proposed fast direct solver also demonstrates its efficiency on multiple right-hand side problems.
Furthermore, the fast direct solver based on the Burton--Miller method outperformed that based on the PMCHWT formulation by about $20\%$ in terms of the computational time.
The proposed method provides a versatile, fast solver, whose performance is independent of the shape of inclusions and computational parameters such as the density, for elastodynamic transmission problems.

As future work, it is important to develop a method that can solve three-dimensional elastic scattering problems in $O(N)$ time.
This is significantly more challenging than the two-dimensional case,
and cannot be achieved by a straightforward extension.
Therefore, the application to three-dimensional problems is the main limitation of the method described here.
The second limitation is the deterioration of accuracy of the fast direct solver
observed in Figures \ref{fig:error_pmchwt_fds} and \ref{fig:error_bm_fds}.
This might not be an issue in practical applications,
but it needs to be resolved from the perspective of solver development.
It also cannot be used for materials that cannot be assumed to be isotropic or elastic.
Furthermore, to analyze composite materials, the proposed method must be extended to handle multiple inclusions.
It is also important to develop low-rank approximation methods with low computational load in two and three dimensions.
One idea is to improve the efficiency of low-rank approximation by taking advantage of the isotropy of the elastodynamic layer potentials \cite{YasuhiroMatsumoto202410-241213}.

\bmhead{Acknowledgements}
This work was supported by
Japan Society for the Promotion of Science under KAKENHI grant number 24K20783
and by the computational resources of TSUBAME4.0 at the Institute of Science Tokyo
provided through the projects 
``Joint Usage/Research Center for Interdisciplinary Large-scale Information Infrastructures (JHPCN)''
and ``High Performance Computing Infrastructure (HPCI)'' in Japan (project IDs jh250045 and jh250070, respectively).

\begin{appendices}

\section{Interpolative factorization} \label{ap:id}
Here, we describe the interpolative factorization.
Let $B^H$ stand for the complex conjugate and transpose of a matrix $B$.
To obtain a left coefficient $U \in \mathbb{C}^{n \times k}$
of a matrix $B \in \mathbb{C}^{n \times m}$, so that $B^{H} \in \mathbb{C}^{m \times n}$,
the interpolative factorization can be used,
where $m, n, k \in \mathbb{N}$; $m \ge n$; and $k \ll m, n$.
Let $B$ be a low-rank approximable matrix.
Note that the interpolative factorization is practically executed
by column-pivoted QR factorization \cite{cheng2005compression}.
For matrix $B$,
the interpolative factorization that computes $U$
based on the column-pivoted QR factorization can be formulated as follows:
\begin{align}
  B^{H} P &= QR \\
  &= \mqty(Q_{11} & Q_{12} \\ Q_{21} & Q_{22}) \mqty(R_{11} & R_{12} \\ & R_{22}) \\
  &\approx \mqty(Q_{11} \\ Q_{21} ) \mqty(R_{11} & R_{12} ) \quad (\text{if } \norm{R_{22}} \text{ is small}) \\
  &= \mqty(Q_{11} R_{11} \\ Q_{21} R_{11}) \mqty(I & R_{11}^{-1} R_{12} ) \\
  B^{H} &= \mqty(Q_{11} R_{11} \\ Q_{21} R_{11}) \mqty(I & R_{11}^{-1} R_{12} )P^{T} \label{eq:brs}
\end{align}
where an orthogonal matrix $Q \in \mathbb{C}^{m \times n}$,
an upper triangular matrix $R \in \mathbb{C}^{n \times n}$,
and a permutation matrix $P$ of size $n \times n$
are obtained from the column-pivoted QR factorization.
In the above relations, we assume that the norm of $R_{22} \in \mathbb{C}^{k \times k}$ is sufficiently small
so that we can low-rank approximate $B$.
The matrices $Q_{11}$, $Q_{12}$, $Q_{21}$, $Q_{22}$, $R_{11}$, and $R_{12}$ are block submatrices of $Q$ and $R$,
whose size is determined so that matrix multiplications can be defined.
Then, to take a complex conjugate and a transpose for both sides of \eqref{eq:brs}, we have
\begin{equation}
  B = \qty{\mqty(I & R_{11}^{-1} R_{12} )P^{T}}^H \mqty(Q_{11} R_{11} \\ Q_{21} R_{11})^H.
\end{equation}
We set $U := \qty{\mqty(I & R_{11}^{-1} R_{12} )P^{T}}^H$, which is the left coefficient of the interpolative factorization for $B$.
Then we have
\begin{equation}
  B = U B_{\mathrm{RS}}, \quad B_{\mathrm{RS}} := \mqty(Q_{11} R_{11} \\ Q_{21} R_{11})^H,
\end{equation}
where $B_{\mathrm{RS}}$ comprises the row skeletons of $B$.
Similarly for $m > n \gg k$,
the right coefficient $V \in {k \times n}$
of the interpolative factorization for $B \in \mathbb{C}^{m \times n}$ can be defined as follows:
\begin{align}
  BP &= QR \\
  &= \mqty(Q_{11} & Q_{12} \\ Q_{21} & Q_{22}) \mqty(R_{11} & R_{12} \\ & R_{22}) \\
  &\approx \mqty(Q_{11} \\ Q_{21} ) \mqty(R_{11} & R_{12} ) \quad (\text{if } \norm{R_{22}} \text{ is small}) \\
  &= \mqty(Q_{11} R_{11} \\ Q_{21} R_{11}) \mqty(I & R_{11}^{-1} R_{12} ) \\
  B &= \mqty(Q_{11} R_{11} \\ Q_{21} R_{11}) \mqty(I & R_{11}^{-1} R_{12} )P^{T} \\
  &= B_{\mathrm{CS}} V, \quad B_{\mathrm{CS}} := \mqty(Q_{11} R_{11} \\ Q_{21} R_{11}), \quad V := \mqty(I & R_{11}^{-1} R_{12} )P^{T},
\end{align}
where the symbols $Q$, $R$, and $P$ are reused, but they actually have different meanings.
In the above relations, $B_{\mathrm{CS}}$ comprises the column skeletons of $B$.

\section{Verification of implementation of the conventional boundary element method} \label{ap:analytical}
The implementation of the conventional boundary element method can be verified by comparing it to the analytical solution on a circle.
The analytical solutions for the displacements $u_1(x)$ and $u_2(x)$ on a circle for transmission problems for elastic waves by a longitudinal plane wave in the two dimensions can be constructed \cite{pao1973Diffraction}.
For the formulations of the boundary element method, the PMCHWT formulation (\eqref{eq:gal_pm1}, \eqref{eq:gal_pm2}),
and the Burton--Miller formulation (\eqref{eq:gal_bm1}, \eqref{eq:gal_bm2}) are employed.

The parameters used in this verification are as follows.
The boundary of an inclusion is a unit circle.
This unit circle is approximated as a polygon with 500 vertices.
This results in both formulations of the boundary element method having 2000 DOFs.
The system of linear equations is solved by the {\it partialPivLU} function in Eigen3 \cite{eigenweb}, which is a standard partial pivoted LU factorization.
The frequency of the incident wave is $\omega = 5$.
The wave speeds in $\Omega_{0}$ are $c_{L}^{(0)} = \sqrt{3}$ and $c_{T}^{(0)} = 1$.
The wave speeds in $\Omega_{1}$ are $c_{L}^{(1)} = 3$ and $c_{T}^{(1)} = 1.5$.
The densities of $\Omega_{0}$ and $\Omega_{1}$ are $\rho^{(0)} = 1$ and $\rho^{(1)} = 2$, respectively.

\begin{figure}[tb]
  \centering
  \includegraphics[width=1.0\linewidth]{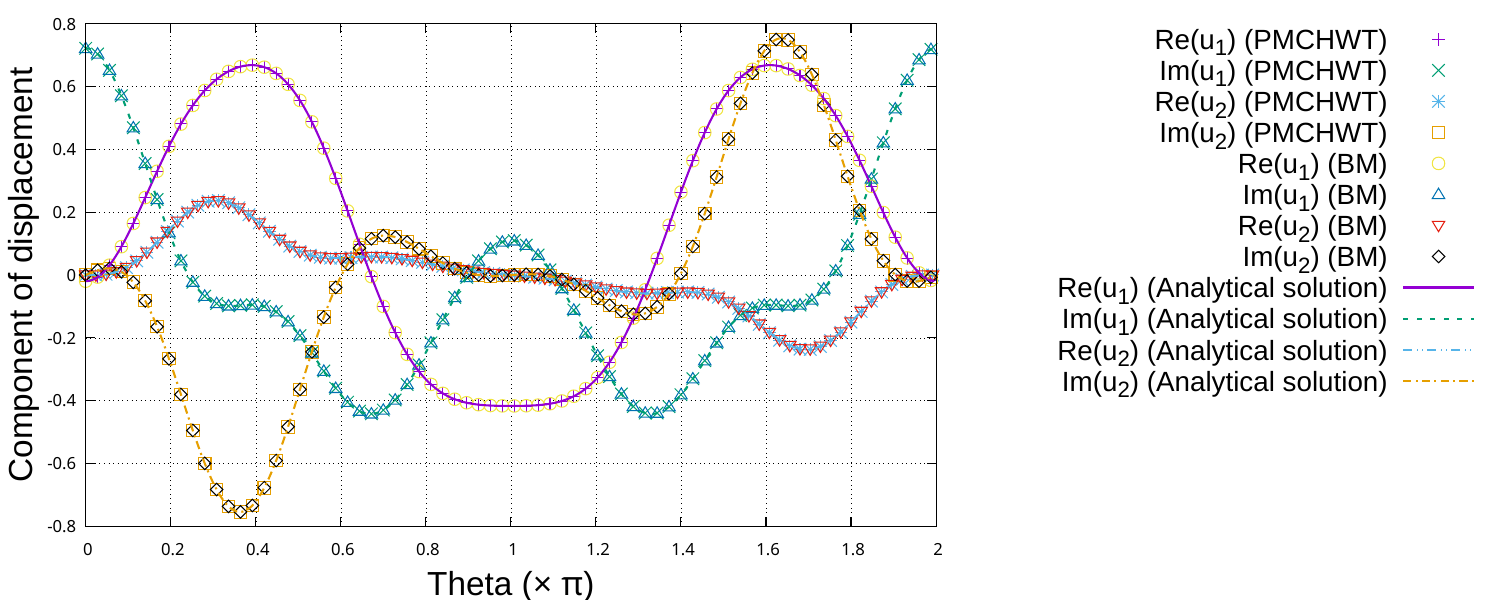}
  \caption{Verifying the implementation of the boundary element method based on the PMCHWT and Burton--Miller (BM) formulation}
  \label{fig:verify_bem}
\end{figure}

Figure \ref{fig:verify_bem} demonstrates good agreement between the computed displacement from both boundary element method formulations and the analytical solution.
In this figure,
the coordinates on the unit circle where each displacement was calculated are expressed as angles in polar coordinates.

\begin{figure}[tb]
  \begin{center}
  \begin{minipage}[b]{0.45\linewidth}
    \centering
    \includegraphics[width=1\linewidth]{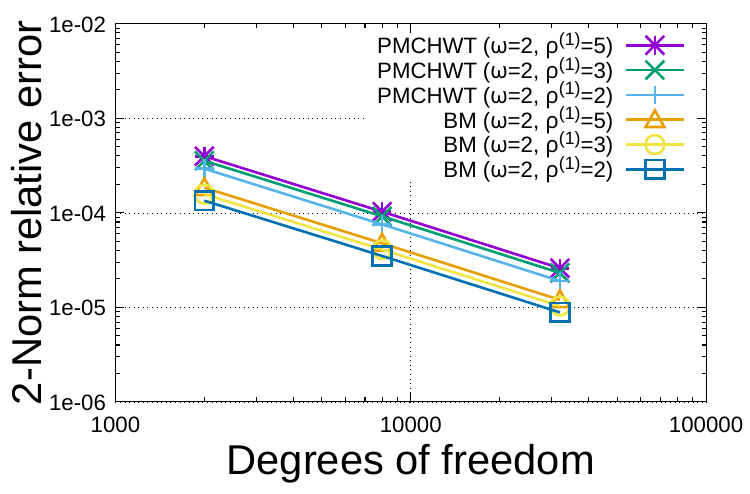}
    \caption{\textcolor{red}{
Two-norm relative error of conventional boundary element method to analytical solution with $\omega = 2$}
}
    \label{fig:verify_o2}
  \end{minipage}
  \end{center}
  \hfill
  \begin{minipage}[b]{0.45\linewidth}
    \centering
    \includegraphics[width=1\linewidth]{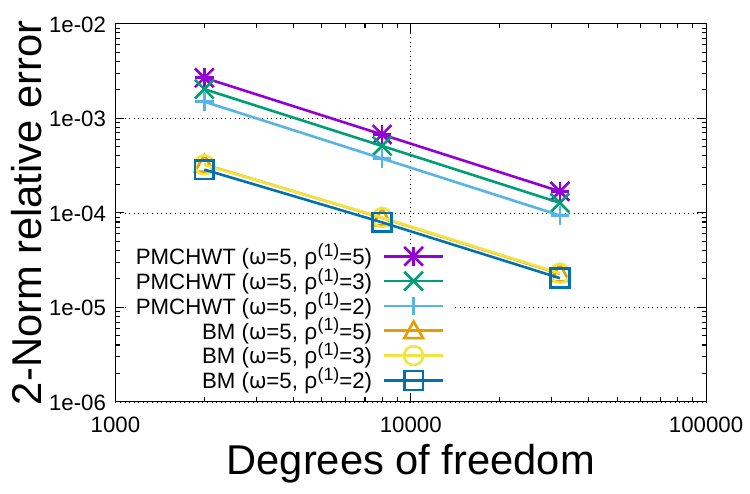}
    \caption{\textcolor{red}{Two-norm relative error of conventional boundary element method to analytical solution with $\omega = 5$}}
    \label{fig:verify_o5}
  \end{minipage}
  \hfill
  \begin{minipage}[b]{0.45\linewidth}
    \centering
    \includegraphics[width=1\linewidth]{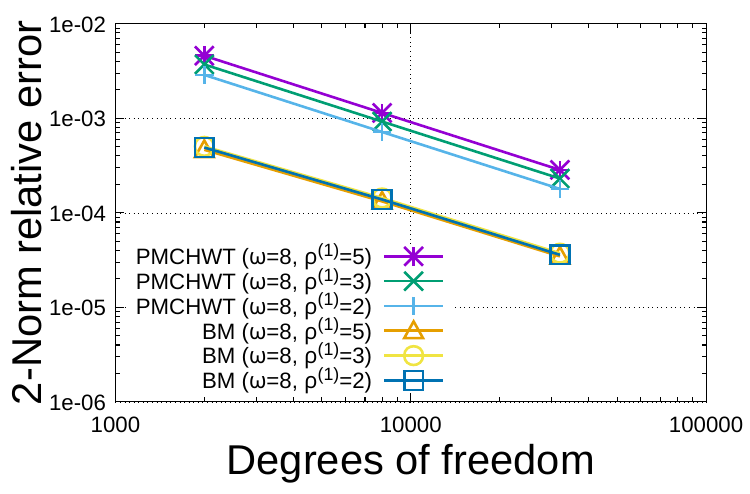}
     \caption{\textcolor{red}{Two-norm relative error of conventional boundary element method to analytical solution with $\omega = 8$}}
    \label{fig:verify_o8}
  \end{minipage}
  \hfill
\end{figure}

Furthermore, we investigate the accuracy of the conventional boundary element method by comparing results across several parameters.
The unit circle is approximated as a polygon with $N$ vertices, where $N = 500, 2000$, and $8000$.
The numerical experiments are performed with frequencies $\omega = 2, 5, 8$ and scatterer densities $\rho^{(1)} = 2, 3, 5$.
Other parameters remain unchanged.
The displacement vector $u^\mathrm{Conv}$ obtained from the conventional boundary element method is expressed as
\[
u^\mathrm{Conv} := ( u_1^1, u_1^2, \ldots, u_1^{N}, u_2^1, u_2^2, \ldots, u_2^{N}) \in \mathbb{C}^{2N},
\]
where $u_i^j$ stands for the $x_i$-directional component of the displacement at the $j$-th vertex of the polygon.
Figures \ref{fig:verify_o2}--\ref{fig:verify_o8} show the relative 2-norm error with respect to the analytical solution, defined by 
\[
  \frac{\norm{u^{\mathrm{Conv}} - {u}^{\mathrm{Analytical}}}_2}{\norm{{u}^{\mathrm{Analytical}}}_2},
\]
where ${u}^{\mathrm{Analytical}}$ is the analytical solution vector on the unit circle calculated at the points corresponding to  those in $u^\mathrm{Conv}$.
These figures confirm the correct implementation of the conventional boundary element method using both PMCHWT and Burton--Miller formulations.

\end{appendices}

\end{document}